\documentclass[12pt, reqno]{amsart}
\usepackage{amssymb,latexsym,amsmath,amscd,amsfonts}
\usepackage{latexsym}

\usepackage{graphicx} 
\usepackage[utf8]{inputenc}
\usepackage{amsmath}
\usepackage{amsthm}
\usepackage{amsfonts}
\usepackage{MnSymbol}
\usepackage{blindtext,amsmath,amsthm,amsfonts, textcomp, mathrsfs, mathtools}
\usepackage{hyperref}
\usepackage{cleveref}
\usepackage[shortlabels]{enumitem}
\numberwithin{equation}{section}

\usepackage[top=30mm,bottom=30mm,left=30mm,right=30mm]{geometry}
\newtheorem{theorem}{Theorem}[section]

\newtheorem{corollary}[theorem]{Corollary}

\newtheorem*{theorem*}{Theorem}
\newtheorem{remark}[theorem]{Remark}
\newtheorem*{remark*}{Remark}
\newtheorem*{problem*}{Problem}
\newtheorem*{conjecture*}{Conjecture}
\newtheorem*{question*}{Question}
\newtheorem{lemma}[theorem]{Lemma}

\newtheorem{definition}[theorem]{Definition}
\newtheorem{example}[theorem]{Example}
\newcommand{\rom}[1]{\uppercase\expandafter{\romannumeral #1\relax}}

\begin{document}

\title{Peripheral Poisson boundary under  compositions and convex combinations}



\author[Bhat and Dias]{B V Rajarama Bhat and Astrid Swizell Dias}



\address[Statistics and Mathematics Unit, Indian Statistical Institute, R V College Post, Bangalore - 560059, India]{}
\email{bvrajaramabhat@gmail.com}
\address[Institute of Mathematical Sciences, Taramani, Chennai]{}
\email{astriddias97@gmail.com}

\begin{abstract}

This article explores the relationship between the peripheral Poisson boundaries of compositions and convex combinations of unital completely positive (UCP) maps on von Neumann algebras and those of their constituent maps. For commuting UCP maps in finite dimensions, the boundary of their compositions is simply the intersection of the  boundaries. Making use of Choquet-Deny theorem it can be seen that, even in infinite dimensions,  any peripheral eigenvector of a convex combination of UCP maps is an eigenvector of all constituent maps with exactly the same eigenvalue. However, the peripheral Poisson boundary may not be a subalgebra, as the extended Choi-Effros products may differ.

Surprisingly, proper convex combinations of peripherally automorphic maps have a boundary entirely independent of the convexity parameters. This class includes UCP maps with a common faithful stationary state and, in particular, mixed unitary channels. Using free groups and the Grigorchuk group, we highlight the obstacles in extending these results to infinite dimensions. Finally, we extend several of these structural results to  continuous time quantum Markov semigroups on full matrix algebras. Finite dimensionality allows us to  provide a concrete description of the peripheral Poisson boundary.

\end{abstract}

\maketitle
\section{Introduction}
Consider a normal, unital completely positive (UCP) map $\tau$ acting on a von Neumann algebra $\mathcal{A}$. The space ${F}(\tau)$ of fixed points, becomes a von Neumann algebra under a modified product, called the Choi-Effros product. The von Neumann algebra ${F}(\tau)$ is called the 
{\em non-commutative Poisson boundary \/} of $\tau$ (See Izumi \cite{Izumi2002} and \cite{Izumi2004}). Similarly, the norm closure of the linear span of peripheral eigenvectors  (this linear span is denoted by $E(\tau )$)  can have an  extended Choi-Effros product,  making it a $C^*$-algebra. This $C^*$-algebra is called the peripheral Poisson boundary of $\tau$. In general, we do not obtain a von Neumann algebra. The notion was introduced by Bhat, Kar and Talwar in \cite{BKT1} and has been studied further  in \cite{BKT2}. It is a very important ingredient to understand the dynamics of the UCP map, as $\tau $ acts as an automorphism of the peripheral Poisson boundary. Further extensions of the theory to contractive completely positive maps and quantum dynamical semigroups can be seen in \cite{BA1}. 
This article investigates how the peripheral Poisson boundaries of convex combinations and compositions of unital completely positive (UCP) maps relate to the boundaries of their constituent maps.

In Section 2, we recall some basic background material. Then in Section 3, to study compositions of UCP maps, we make use of the famous Jacobs–de Leeuw– Glicksberg (JdLG) decomposition. It allows us
to decompose the  von Neumann algebra on which our CP maps are acting as a direct sum of two 
invariant components: a `reversible/almost periodic' part and a `stable/flight' part. The reversible part is seen to be the peripheral Poisson boundary.  We establish that for compositions of commuting UCP maps on finite-dimensional 
$C^*$-algebras, the peripheral Poisson boundary is precisely the intersection of the constituent boundaries (Theorem \ref{pcm} and Corollary \ref{induction}). We also construct counterexamples demonstrating that this  can fail in the absence of either commutativity(Example \ref{noncommuting})
or finite dimensionality (Example \ref{infinite}).

Next we consider convex combinations of UCP maps. Here we are able to get results in two specific contexts: peripherally automorphic maps or mutually commuting maps.  Strong structural results emerge when the constituent maps are peripherally automorphic—that is, when the Choi–Effros product coincides with the original product. In this regime, any convex combination with strictly positive coefficients remains peripherally automorphic. Furthermore, its peripheral Poisson boundary is entirely independent of the chosen convexity parameters and is generated by the linear span of the common peripheral eigenvectors with common eigenvalues (Corollary \ref{bell} and Remark \ref{woods}).

In the broader context of convex combinations of general UCP maps, we are forced to assume that the maps mutually commute. Now we make use of
the well-known Choquet-Deny theorem (See Theorem \ref{choquetdeny}) applicable to abelian semigroups. It shows that under certain situations all harmonic functions are constants. We extend it suitably to the peripheral case in
Theorem \ref{choquetdenyperipheral}. For convex combinations of mutually commuting UCP maps any peripheral eigenvector 
is shown to be a common peripheral eigenvector of the constituents with the same eigenvalue (Corollary \ref{main}). However, surprisingly, the associated peripheral Poisson boundary may fail to be a subalgebra, as the extended Choi–Effros product can differ from that of the constituents. Showing this required a rather subtle construction  (See Example \ref{subtle}).



The notion of peripheral Poisson boundaries was introduced in  \cite{BKT1}. There it was noted that for a normal unital completely positive map $\tau $ acting on a von Neumann algebra, the peripheral Poisson boundary of  $\tau ^n$ does not depend upon $n$ for $n>1$.  Though some basic theory of peripheral Poisson boundary for continuous time quantum Markov semigroups on von Neumann algebras was
developed in \cite{BA1}, this type of time independence of peripheral Poisson boundary  was not shown. Here we are able to do it in finite dimensions with help from JdLG theory (Theorem \ref{ppbQ}). This helps us to extend some results on the peripheral Poisson boundary of compositions of commuting maps to the continuous time, where we are taking sums of commuting generators (Theorem \ref{continuouscomposition}). Unlike the discrete case, time independence  result does not hold in general in infinite dimension due to intricate topological complications (Example \ref{continuouscounter}).


A concrete description of the peripheral Poisson boundary for discrete quantum dynamical semigroups  on the algebra $\mathbb{M}_d$ of complex matrices was given by Amato et. al  \cite{DA}.
In the  last Section we derive the same in continuous time.

\section{Preliminaries}
To begin with, we review the relevant definitions and constructions needed to develop notation for the theory of the peripheral Poisson boundary of unital completely positive maps.

For $d\in \mathbb{N},$ $\mathbb{M}_d$ will denote the algebra of $d\times d$ complex matrices, considered as the space of operators on the Hilbert space $\mathbb{C}^d,$ in the natural way.
Our Hilbert spaces are all complex and separable,
whose inner products, denoted by $\langle \cdot , \cdot \rangle $, are anti-linear
in the first variable. By $\mathscr{B}(\mathcal{H})$ we denote the
algebra of all bounded operators acting on the Hilbert space $\mathcal{H}$. The identity of the Hilbert space $\mathcal{H}$ is denoted by $1_{\mathcal{H}}$. Similarly, we may denote the identity of an algebra $\mathcal{A}$ by $1_{\mathcal{A}}.$ If the Hilbert space or the algebra under consideration is clear from the context, we may drop the subscript. Moreover, elements in $\mathbb{M}_d$, shall be denoted by capital letters, and for those elements lying in $C^*$-algebras other than $\mathbb{M}_d$ or are not specified to be $\mathbb{M}_d$, would be denoted by lowercase letters. In the following, $\mathbb{Z}_+$ and $\mathbb{R}_+$ will denote semigroups of non-negative integers and non-negative real numbers respectively.

We begin with recalling two standard results from the theory of completely positive maps just to set up our notation.

\begin{theorem}\cite{VP}(Stinespring's dilation theorem)\label{Stinespring}
Let $\mathcal{A}$ be a unital $C^*$-algebra, and let $\tau:\mathcal{A}\to\mathscr{B}(\mathcal{H})$ be a completely positive map. Then, there exists a Hilbert space $\mathcal{G}$, a unital $*$-representation $\pi:\mathcal{A}\to\mathscr{B}(\mathcal{G})$, and a bounded operator $V:\mathcal{H}\to\mathcal{G}$, such that $\tau(x)=V^*\pi(x)V, ~\forall x\in\mathcal{A}$, with $||\tau(1_\mathcal{A})||=||V||^2$.
\end{theorem}

We call a triple $(\mathcal{G}, \pi , V)$ as above as a Stinespring dilation of $\tau$. It can be chosen to be {\em minimal}, that is,
$\mathcal{G}=\overline{\pi(\mathcal{A})(\mathcal{H})}$. Minimal Stinespring dilations are unique up to unitary equivalence. It is also to be noted that if $\mathcal{A}$ is a von Neumann algebra and $\tau $ is normal, then the $*$-homomorphsm $\pi $ appearing in its minimal Stinspring dilation is also normal.

Analogous to the notions of absolute continuity and equivalence for measures, there is the following notion for CP maps.
\begin{definition}\label{domination}
Let $\mathcal{A}$ be a $C^*$-algebra, and let $\tau,\eta:\mathcal{A}\to\mathcal{A}$ be completely positive maps. Then,  (i)   $\eta$ is said to be {\em dominated\/} by $\tau$, if $\tau-\eta$ is completely positive; (ii)  $\eta$ is said to be {\em absolutely continuous\/}   to $\tau$, if there exists a scalar $a>0$, such that $a\eta$ is dominated by $\tau$. (iii)  $\eta $ is said to be {\em measure equivalent\/} to $\tau$, if $\eta$
is absolutely continuous to $\tau $ and $\tau $ is absolutely continuous to $\eta .$ \end{definition}

Here is another well-known result on domination of completely positive maps.
\begin{theorem}\cite{WA}(Radon-Nikodym type theorem)\label{RNTT}
Let $\mathcal{A}\subseteq\mathscr{B}(\mathcal{H})$ be a von Neumann algebra, and $\tau:\mathcal{A}\to\mathcal{A}$ be a normal completely positive map, with minimal Stinespring representation $(\mathcal{G},\pi,V)$. Then, a completely positive map $\eta:\mathcal{A}\to\mathcal{A}$ is dominated by $\tau$, if and only if there exists a $D\in\mathscr{B}(\mathcal{G})$ with the following properties:
\item[(1)] $0\leq D\leq 1_\mathcal{G}$,
\item[(2)] $D\pi(x)=\pi(x)D, ~\forall x\in\mathcal{A}$,
\item[(3)] $\eta(x)=V^*\pi(x)DV, ~\forall x\in\mathcal{A}$.
\end{theorem} It follows that if a completely positive map is dominated by a normal completely positive map then it is automatically normal.

Now we recall some basic results on the theory on peripheral Poisson boundary for discrete (refs \cite{BKT1, BKT2}) and continuous time (ref. \cite{BA1}) semigroups of unital
completely positive maps acting on a von Neumann algebra. 
For the discrete-time situation, suppose $\mathcal{A}\subseteq \mathscr{B}(\mathcal{H})$ is a von
Neumann algebra,  and $\tau: \mathcal{A}\to \mathcal{A}$ is a normal,
unital completely positive map.   Consider the set of fixed points of $\tau $:
\begin{equation}\label{Fix1}
F(\tau):=\{x\in \mathcal{A}: \tau(x)=x \}.
\end{equation}
$F(\tau )$ is an operator system, that is, it is a subspace closed under adjoints and contains the identity.  In general, it is not closed under taking products. However, it admits a new product called the Choi-Effros product and this makes $F(\tau )$ a von Neumann algebra. This is known as the  Poisson boundary of $\tau.$ Clearly, $F(\tau )$ is the eigenspace of $\tau$ with respect to eigenvalue 1. Taking a cue from this, for $\lambda \in \mathbb{T}:=\{z: z\in\mathbb{C}, |z|=1\}$, consider
\begin{equation}\label{Fix2}E_\lambda (\tau)=\{x \in \mathcal{A}: \tau(x)=\lambda x \}.\end{equation}  Hence $E_1(\tau )=F(\tau ).$ The elements of $E_{\lambda }$ with $\lambda \in\mathbb{T}$ are known as  peripheral eigenvectors of $\tau.$ 
Define
\begin{equation}\label{Fix3}
E(\tau)=~\mbox{span}\{x\in \mathcal{A}: x\in E_\lambda(\tau) ~\mbox{ for  some } \lambda \in \mathbb{T}\}.
\end{equation}

In studying `quantum dynamics',  in discrete time, one is interested in analyzing the behavior of the discrete semigroup $\{ \tau ^n: n\in \mathbb{Z}_+\}.$ In general $F(\tau ^n)$ varies with $n$. However, an elementary application of the spectral mapping theorem shows (see \cite{BKT1}):
\begin{equation}
E(\tau ^n)= E(\tau ), ~\forall n\geq 1.
\end{equation}
Note that $\tau $ being continuous, $F(\tau )$ is closed. However, this need not be the case with $E(\tau ).$ Take
\begin{equation}\label{Fix4}\mathcal{P}(\tau ) =\overline{E(\tau )}^{\|\cdot \|},\end{equation}
the norm closure of $E(\tau ).$ 

Following the dilation theory of quantum dynamical semigroups as developed in ref. \cite{BhC}  and related papers, in the context of normal, unital completely positive maps on von Neumann algebras, for a triple $(\mathcal{H}, \mathcal{A}, \tau )$, there exists a triple $(\mathcal{K}, \mathcal{B}, \theta )$ called  the {\em minimal dilation} of $\tau$, where (i) $\mathcal{K}$ is a Hilbert space containing $\mathcal{H}$, (ii) $\mathcal{B}\subseteq \mathscr{B}(\mathcal{K})$ is a von Neumann algebra  satisfying $\mathcal{A}= p\mathcal{B}p$, where $p$ is the projection of $\mathcal{K}$ onto $\mathcal{H}$,  and (iii) $\theta : \mathcal{B}\to \mathcal{B}$ is a normal, unital $*$-endomorphism `dilating' $\tau$:
$$\tau ^n(x)=p\theta ^n(x)p, ~~\forall n\geq 0, x\in \mathcal{A}.$$

Here and elsewhere in this dilation theory, the operators $x\in \mathscr{B}(\mathcal{H})$ are identified with $pxp$ in $\mathscr{B}(\mathcal{K})$ for $\mathcal{H}\subseteq \mathcal{K}$. 

As $\theta $ is also a completely positive map in its own right, we already have $E_{\lambda }(\theta ), E(\theta ),$ and $ \mathcal{P}(\theta )$ well defined by simply replacing $\tau $ by $\theta.$ One of the major observations of ref. \cite{BKT1} is that every peripheral eigenvector of $\tau$ `lifts' uniquely to a peripheral eigenvector of $\theta$, and this lift can be extended to the
whole of $\mathcal{P}(\tau )$.  More specifically,  given $x\in \mathcal{P}(\tau )$, there exists a unique $\widehat{x}$ in $\mathcal{P}(\theta )$ such that $x=p\widehat{x}p.$ 
The mapping $x\mapsto \widehat{x}$ from $\mathcal{P}(\tau )$ to $\mathcal{P}(\theta )$ is a bijection. In fact,   it is an isometric isomorphism of the operator systems. The inverse map of this is nothing but the compression map $z\mapsto pzp$, restricted to $\mathcal{P}(\theta ).$  Now, with $\theta $ being a unital $*$-endomorphism, $\mathcal{P}(\theta )$ is seen to be a unital $C^*$-subalgebra of $\mathcal{B}$.  This allows us to define a product, denoted by `$\circ $', on $\mathcal{P}(\tau )$ by taking
$x\circ y = p(\widehat{x}.\widehat{y})p.$

Inducing the extended Choi-Effros product on $\mathcal{P}(\tau )$, we have $(\mathcal{P}(\tau ), \circ )$ as a unital $C^*$-algebra and is called as the peripheral Poisson boundary of $\tau $. In fact, for $x\in E_{\lambda }(\tau ), ~ y\in E_{\mu } (\tau )$, where $\lambda,\mu \in \mathbb{T}$, 
we have the formula: \begin{equation}\label{limit formula} x\circ y = ~\mbox{s.}\lim _{n\to \infty}\frac{\tau ^n (xy)}{\lambda ^n\mu ^n }\end{equation} where `$\mbox{s.}$' denotes the strong operator topology.

Replacing the semigroup $\mathbb{Z}_+$ by the semigroup $\mathbb{R}_+$, we have the continuous-time semigroup of contractive CP maps called quantum dynamical semigroups.

\begin{definition}
A quantum 
dynamical semigroup is a family of linear maps $\tau =\{\tau _t: t\in
\mathbb{R}_+\}$ acting on a von Neumann algebra $\mathcal{A}\subseteq \mathscr{B}(\mathcal{H})$, where $\mathcal{H}$ is a complex, separable
Hilbert space, 
satisfying:
\begin{itemize}
\item[1)] $\tau_t$ is completely positive, $\forall t\in \mathbb{R}_+$;
\item[2)] $\tau_s(\tau_t(x))=\tau_{s+t}(x), ~~\forall s,t\in \mathbb{R}_+, x\in
\mathcal{A}$;
\item[3)] $\tau_t(1) \leq 1$, $\forall t\in \mathbb{R}_+$;
\item[4)] $x \mapsto \tau_t(x)$ is normal, $\forall t\in
\mathbb{R}_+$;
\item[5)]$t \mapsto \tau_t(x)$ is continuous in strong operator topology for every
$x\in \mathcal{A}.$
\end{itemize}\end{definition}

A quantum dynamical semigroup is called {\em unital/Markov\/}, when $\tau_t(1)=1$, $\forall t$. It is called an $E$-semigroup, if $\tau _t$ is a $*$-endomorphism for every $t$ and an $E$-semigroup is called an $E_0$-semigroup, if it is unital.

The concept of fixed points and peripheral eigenvectors in the continuous case, is analogous to the discrete case. Consider a quantum dynamical semigroup $\tau $ acting on a von Neumann algebra $\mathcal{A}$. 

{(i)} The set of harmonic elements or fixed points  of $\tau$, is given by  
\begin{equation}\label{fixed}
F(\tau):=\{x\in\mathcal{A}: \tau_t(x)=x;~ \forall t \in \mathbb{R}_+\}. \end{equation}

{(ii)} The {\em peripheral eigenspace\/} of $\tau $ is given by \begin{equation}\label{Peripheral} E_{a}(\tau):=\{x\in \mathcal{A}: \tau_t(x)=\exp{(iat)} x;~ \forall t \in \mathbb{R}_+ \}, \end{equation}  with $a \in \mathbb{R}$.  
Just as in the discrete case, we consider 
\begin{equation}\label{Fix6} E(\tau):=  \mbox{span}\{x\in E_a( \tau) :  a\in \mathbb{R} \}.\end{equation}
\begin{equation}\label{Fix7}
\mathcal{P}(\tau):=\overline{\mbox{span}}{\{x\in E(\tau)\}}.\end{equation}

Analogously, one has the extended Choi-Effros product on $\mathcal{P}(\tau )$, and again $(\mathcal{P}(\tau ), \circ )$ is a unital $C^*$-algebra and is called as the peripheral Poisson boundary of $\tau $.  For $x\in E_{a }(\tau ), ~ y\in E_{b } (\tau )$, where $a,b \in \mathbb{R}$, 
we have: \begin{equation}\label{limit formula 2} x\circ y = ~\mbox{s.}\lim _{t\to \infty}\frac{\tau _t (xy)}{\exp{i(a+b)t} }.\end{equation}



%


In the following we have multiple occasions where we need to consider semigroups of completely positive maps which are not necessarily unital.
For instance, we need to look at adjoints of quantum Markov semigroups, which are of course trace preserving but not necessarily unital. For such semigroups  also we continue to use the notation introduced from (\ref{Fix1}) to (\ref{Fix7}).
In some such cases,  $\mathcal{P}(\tau)$  may not have any  $C^*$-algebraic structure on it, but it is convenient to retain the same notation.

Next, we dive into key definitions regarding the operator semigroups and state the fundamental Jacobs–de Leeuw–Glicksberg (JdLG) decomposition theorem (ref. \cite{Emel}). The  weak topology appears naturally in this theory, and is mentioned below to be faithful to the original literature. But in the present article the concept is applied only for finite dimensional spaces. Hence we may ignore such topological subtleties.

\begin{definition}
For a given semigroup $J$, a multiplicative homomorphism $\alpha: J\to \mathbb{T}$ is called a character. 
\end{definition}
\begin{definition} Let the semigroup $J=\mathbb{R}_+ ~\mbox{or}~ \mathbb{Z}_+$,
and $T=\{T_t: t\in J\}$ be a  semigroup of operators on the Banach space $\mathcal{X}.$  Then, $T$ is called (weakly) almost periodic, if the orbit $\{T_t x: t\in J \}$ is (weakly) pre-compact for every $x\in \mathcal{X}.$ 
\end{definition}

Consider the semigroup $J=\mathbb{R}_+ ~\mbox{or}~ \mathbb{Z}_+$ and the semigroup of operators $T=\{T_t:t\in J\}$, acting on the Banach space $\mathcal{X}$. 
Define the two subspaces of $\mathcal{X}$ as follows:
$$\mathcal{X}_{fl}(T)=\{x\in \mathcal{X} \ : \ 0 \in \overline{\{T_t(x):t\in J\}}^w\},$$
$$\mathcal{X}_r(T)=
\{x\in \mathcal{X}: y\in \overline{\{T_t(x): t\in J\}}^w \ \mbox{implies}\ x\in \overline{\{T_t(y): t\in J\}}^w\},$$
where $\overline{\mathcal{Y}}^w$ denotes weak-closure on the subspace $\mathcal{Y}$ of $\mathcal{X}.$
The vectors in $\mathcal{X}_{fl}(T)$ are called {\em flight\/} vectors and the vectors in $\mathcal{X}_r(T)$ are called  {\em reversible\/} vectors (\cite{Emel}). 

Recall that, a one-parameter semigroup $T=\{T_t\}_{t\in \mathbb{R}_+}$ of bounded linear operators on  a Banach space $\mathcal{X}$ is called a ${C}_0$ semigroup, if $T_0=id.$, where $id.$ is the identity map and it is point-wise continuous, that is, 
$t \mapsto T_t (x)$ from $\mathbb{R}_+ $ to $\mathcal{X}$ is continuous for every  $x\in \mathcal{X}$. Moreover, the semigroup $T=\{T_t\}_{t\in J}$ is called bounded, if $M_T=\mbox{sup}\{||T_t||: \ t\in J\}<\infty$.

An interesting result to be acquainted with is as follows.
\begin{lemma}(Theorem 1.1.2 in \cite{Emel})\label{Reflex}
A bounded ${C}_0$ semigroup  $T=\{T_t\}_{t\in \mathbb{R}_+}$ acting on a reflexive  Banach space $\mathcal{X}$ is weakly almost periodic. 
\end{lemma}

We are now prepared to state the Jacobs-de Leeuw-Glicksberg Decomposition theorem.
\begin{theorem}\label{Jacob}
(Jacobs-de Leeuw-Glicksberg (JdLG) Decomposition \cite{Emel}) Given a weakly almost periodic semigroup $T=\{T_t\}_{t\in J}$ ($J=\mathbb{R}_+ ~\mbox{or}~ \mathbb{Z}_+$) on a Banach space $\mathcal{X}$,
 then $\mathcal{X}= \mathcal{X}_r(T)\oplus \mathcal{X}_{fl}(T),$
and $\mathcal{X}_r(T)=\overline{\mbox{span}}\{x\in \mathcal{X}\ : \ \mbox{there exists a character } \alpha: J \to \mathbb{T } \ \mbox{such that} \  \tau_t (x) = \alpha(t)x, \ \forall t\}$. 
If $T$ is almost periodic, then $\mathcal{X}_{fl}(T)= \{ x\in \mathcal{X}: \ \lim _{t\to \infty}||T_t(x)||=0\}$.
\end{theorem}

We can establish an additional structure for the reversible and flight spaces in the above theorem. Before presenting this result, we recall a few definitions.

\begin{definition}\cite{EFHN}
Let $S$ be a non-empty semigroup.
\item[(1)] A right (left, two-sided) ideal of a semigroup $S$ is called minimal, if it is minimal with respect to the set inclusion within the set of all right (left, two-sided) ideals.
\item[(2)] An idempotent element in a minimal right (left) ideal is called a minimal idempotent. 
\end{definition}
Moreover, minimal right (left) ideals are either disjoint or equal.

\begin{definition}\cite{EFHN}
Let $\mathcal{X}$ be a Banach space and $\mathcal{J}\subseteq\mathscr{B}(\mathcal{X})$ be a semigroup of operators with the weak operator topology. Then, $\mathcal{J}$ is called Jacobs-de Leeuw-Glicksberg (JdLG) admissible, if $\mathcal{J}$ is relatively weakly compact and the weak-closure contains a unique minimal idempotent. 
\end{definition}

\begin{remark}\label{JDLG-admissible}
Let $\mathcal{X}$ be a Banach space and $T=\{T_t\}_{t\in J}\subseteq\mathscr{B}(\mathcal{X})$ be a semigroup of operators as in Theorem \ref{Jacob}. Since, $T=\{T_t\}_{t\in J}$ is abelian, by Theorem 16.19 in \cite{EFHN}, $\mathcal{J}$ is JdLG admissible. Then, there exists an unique minimal idempotent $Q\in \overline{\mathcal{J}}^w$, such that $\mathcal{X}_{fl}(\mathcal{J})=Ker(Q)$ and $\mathcal{X}_{r}(\mathcal{J})=Ran(Q)$, with $\mathcal{X}=\mathcal{X}_{fl}(\mathcal{J})\oplus\mathcal{X}_{r}(\mathcal{J})$ (see pages 308,309 in ref. \cite{EFHN}).
\end{remark}

\section{Peripheral Poisson boundary of compositions of UCP maps}\label{COMPCOMB}

In this section, we explore the conditions under which, the peripheral Poisson boundary of $\tau$, which is a compositions of UCP maps, acting on a finite dimensional $C^*$-algebra, is related to the peripheral Poisson boundaries of the UCP maps under composition. 

\begin{remark}\label{rev}
Let $\mathcal{A}$ be a finite dimensional $C^*$-algebra, $\tau:\mathcal{A}\to\mathcal{A}$ be a map such that $\mathcal{J}_\tau:=\{id.,\tau,\tau^2,\tau^3,\ldots,\tau^n,\ldots\}$ be an uniformly bounded semigroup. Then, $\mathcal{J}_\tau$ is an abelian, relatively compact semigroup. By Theorem 16.19 in \cite{EFHN}, $\mathcal{J}_\tau$ is JdLG-admissible. Under these conditions, using Corollary 16.32 in \cite{EFHN}, we have that $\mathcal{A}_r(\mathcal{J}_\tau)=\overline{span}\{x\in\mathcal{A}:\forall n\in \mathbb{N},~ \mbox{there exists a} ~\lambda\in\mathbb{T}, ~\mbox{such that }\tau^n(x)=\lambda ^n x\}$. Hence, by Theorem 4.3 in \cite{BKT1}, $\mathcal{A}_r(\mathcal{J}_\tau)=\mathcal{P}(\tau)$.
\end{remark}

Having the required essential theoretical framework, we observe the following theorem.

\begin{theorem} \label{pcm}
Let $\mathcal{A}$ be a finite dimensional $C^*$-algebra and  $\tau,\alpha,\beta:\mathcal{A}\to\mathcal{A}$ be UCP maps such that $\alpha$ and $\beta$ commute with each other and $\tau=\alpha\beta$. Then, $\mathcal{P}(\tau)= \mathcal{P}(\alpha)\cap \mathcal{P}(\beta)$. 
\end{theorem}
\begin{proof}
Consider the semigroups $\mathcal{J}_\tau=\{id.,\tau,\tau^2,\ldots\}$, $\mathcal{J}_\alpha=\{id.,\alpha,\alpha^2,\ldots\}$ and  $\mathcal{J}_\beta=\{id.,\beta,\beta^2,\ldots\}$ on the Banach space $\mathcal{A}.$ Clearly, they are all abelian. Moreover, due to finite dimensionality of $\mathcal{A}$, and contractivity of $\tau , \alpha , $ and $\beta$, they are relatively compact. Hence, by Theorem 16.19 of \cite{EFHN},  $\mathcal{J}_\tau,\mathcal{J}_\alpha$ and $\mathcal{J}_\beta$ are JdLG-admissible. Let $R,P, \ \mbox{and} \ Q$ be the unique minimal idempotents for $\overline{\mathcal{J}}_\tau$, $\overline{\mathcal{J}}_\alpha$, and $\overline{\mathcal{J}}_\beta$ respectively. 
By Remark \ref{JDLG-admissible} and Remark \ref{rev}, $\mathcal{P}(\tau)=\mathcal{A}_r(\mathcal{J}_\tau)=Ran(R)$, $\mathcal{P}(\alpha)=\mathcal{A}_r(\mathcal{J}_\alpha)=Ran(P)$, and $\mathcal{P}(\beta)=\mathcal{A}_r(\mathcal{J}_\beta)=Ran(Q)$. Since $\tau$, $\alpha$, and $\beta$ mutually commute, $R,P,$ and $Q$ also mutually commute. To see this, note that by  a known result (see Lemma 2.2 in \cite{BKT1}), if  $\{\lambda_1,\lambda_2,\ldots\lambda_p\}$ is a finite subset of $\mathbb{T}$, then there exists a strictly increasing sequence of  natural numbers, say $\{n_k\}_{k\in\mathbb{N}}$, such that  $\lim_{k\to\infty}\lambda_j^{n_k}=1$ for every $1\leq j\leq p.$ Then, if $\{\lambda_1,\lambda_2,\ldots\lambda_p\}$ contains all peripheral eigenvalues of $\tau,\alpha,$ and $\beta$,  we have $\lim_{k\to\infty}\tau^{n_k}=R,~\lim_{k\to\infty}\alpha^{n_k}=P,~\lim_{k\to\infty}\beta^{n_k}=Q$. Hence, $\lim_{k\to\infty}(\tau\alpha)^{n_k}=RP=PR=\lim_{k\to\infty}(\alpha\tau)^{n_k}$, $\lim_{k\to\infty}(\tau\beta)^{n_k}=RQ=QR=\lim_{k\to\infty}(\beta\tau)^{n_k}$ and $\lim_{k\to\infty}(\alpha\beta)^{n_k}= PQ=QP=\lim_{k\to\infty}(\beta\alpha)^{n_k}$. Thereby, $R=\lim_{k\to\infty}\tau^{n_k}=\lim_{k\to\infty}(\alpha\beta)^{n_k}=PQ$. 

Since $R^2=R$ and $PR=R$, we have $\mathcal{P}(\tau)\subseteq\mathcal{P}(\alpha)$. Moreover, since $R^2=R$ and $QR=RQ=R$, implies $\mathcal{P}(\tau)\subseteq\mathcal{P}(\beta)$. Therefore, $\mathcal{P}(\tau)\subseteq\mathcal{P}(\alpha)\cap\mathcal{P}(\beta)$.
If $x\in \mathcal{P}(\alpha)\cap\mathcal{P}(\beta)$, then $Rx=PQx=Px=x$, therefore, $x\in \mathcal{P}(\tau)$, implying $\mathcal{P}(\alpha)\cap\mathcal{P}(\beta)\subseteq \mathcal{P}(\tau)$.
Hence, the claim.
\end{proof}

A simple mathematical induction leads to  the following result.

\begin{corollary}\label{induction}
Let $\mathcal{A}$ be a finite dimensional $C^*$-algebra, and let $\tau,\alpha_1,\alpha_2,\ldots\alpha_m:\mathcal{A}\to\mathcal{A}$ be UCP maps such that $\tau=\alpha_1\alpha_2\cdots\alpha_m$, and all the $\alpha_i$'s mutually commute. Then, $ \mathcal{P}(\tau)= \bigcap_{i=1}^m\mathcal{P}(\alpha_i)$. 
\end{corollary} 

Note that if UCP maps do not mutually commute, the peripheral Poisson boundary of their compositions may not be equal to the intersection of the peripheral Poisson boundaries. We give below an example to illustrate this.

\begin{example}\label{noncommuting}
Let $\tau,\alpha,\beta:\mathbb{C}^3\to\mathbb{C}^3$ be UCP maps defined as  $\tau(x_1,x_2,x_3)=\alpha(x_1,x_2,x_3)=(x_1,x_2,x_2)$ and $\beta(x_1,x_2,x_3)$ $=(x_1,x_2,\frac{1}{3}(x_1+x_2+x_3))$. Clearly $\alpha\beta\neq\beta\alpha $. Note that $\tau=\alpha\beta$. The peripheral Poisson boundaries for $\tau,\alpha,\beta$ are   $\mathcal{P}(\tau)=~\mbox{span}\{(1,0,0),(0,1,1)\}=\mathcal{P}(\alpha)$, and $ \mathcal{P}(\beta)=~\mbox{span}\{(1,1,1),(1,-1,0)\}.$ Clearly, $\mathcal{P}(\alpha)\cap\mathcal{P}(\beta)=\mathbb{C}(1,1,1)\neq \mathcal{P}(\tau)$.
\end{example}

Next, we ask whether the mutual commutativity of the $\alpha_i$'s still suffices for the infinite-dimensional von Neumann algebra, to have the peripheral Poisson boundary of $\tau$ coincide with the intersection of the peripheral Poisson boundaries of the $\alpha_i$'s. This is not the case. We demonstrate this with the example below.

\begin{example}\label{infinite}
Consider the von Neumann algebra $\mathcal{A}=L^\infty(\mathbb{T})$. Let the UCP maps $\alpha_1,\alpha_2:\mathcal{A}\to\mathcal{A}$ be defined as $\alpha_1(f)(z)=f(\exp{(2\pi i\gamma_1)}z)$, and $\alpha_2(f)(z)= f(\exp{(2\pi i\gamma_2)}z)$ for $\gamma_1,\gamma_2\in \mathbb{R}$ such that $\gamma_1,\gamma_2\notin \mathbb{Q}$, but $\gamma_1+\gamma _2\in \mathbb{Q}$. Define the UCP map $\tau:\mathcal{A}\to\mathcal{A}$ as $\tau=\alpha_1\alpha_2$, i.e. $\tau(f)(z)=f(\exp{(2\pi i\gamma)}z), \ \mbox{where} \ \gamma=\gamma_1+\gamma_2$. From Example \ref{continuouscounter} below, we observe that the peripheral Poisson boundary of $\tau$,  $\mathcal{P}(\tau)=L^\infty(\mathbb{T})$,  where as $\mathcal{P}(\alpha_1)=C(\mathbb{T})=\mathcal{P}(\alpha_2)$.
\end{example}

\section{Peripheral Poisson boundary of convex combinations of UCP maps}\label{CONVEXCOMB}

In this section, we are interested in determining the peripheral Poisson boundary of convex combinations of UCP maps.
\begin{example}\label{fif}
Let $A,B \in \mathbb{M}_d \ \mbox{such} \ \mbox{that} \ 0 \leq A \leq B$ and $B$ is unitarily equivalent to $A$. Considering the trace it is easy see that this implies $A=B$.  However, this is not true in general for the infinite dimensional case. Let $\{a_n\}_{n \in \mathbb{Z}}$ be positive real numbers such that $0<a_n < a_{n+1} <1, \ \forall n$. Let $\mathcal{H}$ be a Hilbert space with the orthonormal basis $\{e_n: n \in \mathbb{Z}\}$. Let $U:\mathcal{H}\to\mathcal{H}$ be the bilateral shift defined by $U(e_n)=e_{n+1}, \ \forall n$. Let $A$ be the diagonal operator such that $A(e_n)=a_n e_n, \ \forall n$.  Take $B=U^* A U$. Note that $\langle e_n , B e_n \rangle = \langle U e_n , A U e_n \rangle = a_{n+1}$. So, it is easy to see that B is a diagonal operator such that $B e_n = a_{n+1} e_n, \ \forall n$. Therefore, $A \leq B$  and  $A$ is unitarily equivalent to $B$, but still $A\neq B.$ 
\end{example}

The example above conveys that under unitary conjugation some of the order  properties holding in finite dimensions, do not necessarily follow for the infinite dimensional case. We will see manifestations of this in the context of convex combinations of unitary conjugations.
In $\mathbb{M}_d$, convex combinations of conjugations by unitaries are known as mixed unitary channels. Note that they are unital and trace preserving. The following result is folklore.

\begin{theorem}\label{finite unitary} Let $\tau:\mathbb{M}_d\to \mathbb{M}_d$ be a mixed unitary map, that is, $\tau(X)=\sum_{i=1}^n p_iU_i^*XU_i$, where $U_i$'s are unitaries,  $ 0< p_i < 1, \ \forall i$ and  $\sum_{i=1}^n p_i=1.$ If $X$ is a peripheral eigenvector, that is,  $\tau(X)=\lambda X $ for $\lambda\in \mathbb{T}$, then $U_i^*XU_i=\lambda X, \ \forall i$.
\end{theorem}
\begin{proof}
By assumption, $\tau$ is a stationary map with the faithful invariant state $tr=\frac{1}{d}Tr$, hence, by Corollary 2.8 in \cite{BKT2}, $\tau$ is peripherally automorphic. 
By Theorem 2.6 in \cite{BKT2}, $\tau(X)=\lambda X $ for $\lambda\in \mathbb{T}$, implies $XU_i=\lambda U_i X, \ \forall i$. Therefore, the claim.
\end{proof}
 Here we have a significant generalization of this result applicable for peripherally automorphic maps on general von Neumann algebras. We remind the reader that a UCP map $\tau $ is called peripherally automorphic if the extended Choi-Effros product on its peripheral Poisson boundary is same as the original product.

\begin{theorem}\label{extn}
Let $\mathcal{A}\subseteq\mathscr{B}(\mathcal{H})$ be a von Neumann algebra, and $\tau:\mathcal{A}\to\mathcal{A}$ be a normal UCP map. Consider a minimal Stinespring representation  $(\mathcal{G},\pi,V)$ of $\tau $. 
Then the following statements are equivalent:
\item[(1)] $\tau$ is peripherally automorphic.
\item[(2)] For  every $y\in E_\lambda(\tau),$ with $  \lambda\in\mathbb{T}$, $\tau(y^*y)=y^*y$.
\item[(3)] For $\lambda\in\mathbb{T}$, $y\in E_\lambda(\tau)$ iff $\pi(y)V=\lambda Vy$.
\end{theorem}

\begin{proof} Let (1) be true. For $\lambda\in\mathbb{T}$, if
$y\in E_\lambda(\tau)$, $y^*\circ y=y^*y$. Therefore $\tau(y^*y)=\tau(y^*\circ y)=\tau(y^*)\circ\tau(y)=|\lambda |^2y^*\circ y=y^*y$. Hence, $(2)$ follows from (1).

Now let (2) be true. Let $y\in E_\lambda(\tau), ~\mbox{for}~  \lambda\in\mathbb{T}$. Then, $(\pi(y)V-\lambda Vy)^*(\pi(y)V-\lambda Vy)=V^*\pi(y^*y)V-\lambda V^*\pi(y^*)Vy-\overline{\lambda}y^*V^*\pi(y)V+y^*V^* Vy=\tau(y^*y)-\lambda\tau(y^*)y-\overline{\lambda}y^*\tau(y)+y^*y=0$. Hence
$\pi(y)V=\lambda Vy$. Conversely, let $\pi(y)V=\lambda Vy$, then $\tau(y)=V^*\pi(y)V=\lambda V^*Vy=\lambda y$. Therefore, (2) implies (3).

Assume that (3) holds.
Consider $x\in E_\lambda(\tau),~y\in E_\mu(\tau)$ for $\lambda,\mu\in\mathbb{T}.$ Note that $\tau(xy)=V^*\pi(xy)V=V^*\pi(x)\pi(y)V=\mu V^*\pi(x)Vy=\lambda\mu xy$. To show that $\tau$ is peripherally automorphic, it is required to show that $x\circ y=xy$. This is now easy to
see from the formula (\ref{limit formula}): $x\circ y=\mbox{s.}\lim_{n\to\infty}\frac{\tau^n(xy)}{\lambda^n\mu^n}=xy$. Hence, (1) follows from (3). 
\end{proof}

In the following theorem we use the notions of `absolute continuity' and  `measure equivalence' as introduced in Definition \ref{domination} to compare CP maps.
\begin{theorem}\label{wonder}
Let $\mathcal{A}\subseteq\mathscr{B}(\mathcal{H})$ be a von Neumann algebra and $\tau,\eta:\mathcal{A}\to\mathcal{A}$ be normal UCP maps. Assume $\tau$ is peripherally automorphic and $\eta$ is absolutely continuous  to $\tau$. Then, for $\lambda\in \mathbb{T},$ $\tau(x)=\lambda x $, implies $\eta(x)=\lambda x$ and in particular, $\mathcal{P}(\tau )\subseteq \mathcal{P}(\eta ).$ If further, $\eta $ is measure equivalent to $\tau$, then $\eta$ is peripherally automorphic and $\mathcal{P}(\eta )=\mathcal{P}(\tau ).$ 
\end{theorem}
\begin{proof}
Consider $a>0$, such that $a\eta$ is dominated by $\tau$. Let $(\mathcal{G},\pi,V)$ be  the minimal Stinespring representation of $\tau$. By Theorem \ref{RNTT}, there exists a unique  positive contraction $D\in \mathscr{B}(\mathcal{G}),$  such that $a\eta(x)=V^*\pi(x)DV,$ and $\pi(x)D=D\pi(x),\ \forall x\in\mathcal{A}$. Using Theorem \ref{extn}, we have that if $\tau(x)=\lambda x; \ \lambda\in\mathbb{T}$, since 
$\pi(x)V=\lambda Vx$, $\eta(x)=\frac{1}{a}V^*D\pi(x)V=\lambda(\frac{1}{a}V^*DV)x=\lambda x$. 
Now the second part follows, once we show that $\eta $ is also peripherally automorphic. Let $x\in E_\lambda(\eta), ~y\in E_\mu(\eta), \mbox{ for} ~\lambda,\mu\in\mathbb{T}$. Since, $\eta(xy)=\frac{1}{a}V^*D\pi(xy)V=\frac{1}{a}V^*D\pi(x)\pi(y)V=\lambda\mu (\frac{V^*DV}{a})xy$, $\eta(xy)=\lambda\mu xy$. Therefore, $x\circ_\eta y=xy$, where $\circ_\eta$ is the extended Choi-Effros product on the peripheral Poisson boundary $\mathcal{P}(\eta) $. Since, $x$ and $y$ were arbitrarily chosen, $\eta$ is as claimed.
\end{proof}

Using the theorem stated above, and under the required assumptions, we establish a connection between the peripheral Poisson boundary of a convex combination of normal UCP maps and the peripheral Poisson boundaries of its components.

\begin{corollary} \label{bell}
Let $\mathcal{A}\subseteq\mathscr{B}(\mathcal{H})$ be a von Neumann algebra and $\tau:\mathcal{A}\to\mathcal{A}$ be a normal UCP map,  with a faithful stationary state $\varphi$. Assuming $\tau$ to be a convex combination of the normal UCP maps: $\tau_i$'s, i.e. $\tau(x)=\sum_{i=1}^n p_i\tau_i(x)$, where $\sum_{i=1}^n p_i=1$ and $ 0< p_i < 1, \ \forall i$. Then, for $\lambda\in \mathbb{T}$, $\tau(x)=\lambda x $ iff $\tau_i(x)=\lambda x, \ \forall i$.
Consequently, $E_\lambda(\tau)=\bigcap_{i=1}^n E_\lambda(\tau_i)$ for $\lambda\in\mathbb{T}$ and $\mathcal{P}(\tau )
= ~\overline{\mbox{span}}\bigcup _{\lambda \in \mathbb{T}} \bigcap_{i=1}^n E_\lambda(\tau_i)
\subseteq \bigcap _{i=1}^n\mathcal{P}(\tau_i).$
\end{corollary}
\begin{proof}
Note that $\tau$ is peripherally automorphic. This is because, for $x\in E_\lambda(\tau), ~\lambda\in\mathbb{T}$, by Kadison-Schwarz inequality, $\tau(x^*x)\geq\tau(x^*)\tau(x)$, which implies $\tau(x^*x)\geq x^*x$, and by assumption, $\varphi(\tau(x^*x)- x^*x)=0$ implies $\tau(x^*x)=x^*x$.
For every $1\leq i\leq n$, $p_i\tau_i$ is dominated by $\tau$. By Theorem \ref{wonder}, $\tau(x)=\lambda x $ for $\lambda\in \mathbb{T}$, implies $\tau_i(x)=\lambda x, \ \forall i$. The second part follows easily.
\end{proof}

Observe that the above corollary means that for the peripherally automorphic $\tau$, the peripheral Poisson boundary of $\tau $ is actually independent of the convexity parameters $p_1, p_2, \ldots,p_n$ as long as they are strictly positive (such convex combinations are known as proper convex combinations). 
We ask whether, we can drop the  the peripherally automorphic assumption and arrive at same conclusions.  The answer is no. We demonstrate this in the example below.
 
\begin{example}
 Consider the von Neumann algebra $\mathcal{A}=\mathbb{M}_2\oplus l^\infty (\mathbb{Z}_+)$, and normal UCP maps $\tau_1,\tau_2:\mathcal{A}\to\mathcal{A}$, where 
 $$\tau_1 (\begin{bmatrix}
     x_{00}&x_{01}\\x_{10}&x_{11}
 \end{bmatrix}\oplus (y_0,y_1,y_2,\ldots))=\begin{bmatrix}
     x_{00}&x_{01}\\x_{10}&x_{11}
 \end{bmatrix}\oplus (x_{00},y_0,y_1,y_2,\ldots);$$
 $$\tau_2 (\begin{bmatrix}
     x_{00}&x_{01}\\x_{10}&x_{11}
 \end{bmatrix}\oplus (y_0,y_1,y_2,\ldots))=\begin{bmatrix}
     x_{00}&x_{01}\\x_{10}&x_{11}
 \end{bmatrix}\oplus (x_{11},y_0,y_1,y_2,\ldots).$$
 Define the normal UCP maps $\tau$ and $\eta$ acting on $\mathcal{A}$, which are convex combinations of the UCP maps $\tau_1,\tau_2$ as $\tau=\frac{1}{2}(\tau_1+\tau _2)$ and $\eta=\frac{1}{4}(\tau_1+3\tau_2).$ 
%
%
Here, $1$ happens to be the only peripheral eigenvalue for $\tau$ and $\eta$. The corresponding peripheral Poisson boundaries are $\mathcal{P}(\tau)=span\{E_{01}\oplus 0,E_{10}\oplus 0,(E_{00}-E_{11})\oplus 0, (E_{00}+E_{11})\oplus (1,1,1,\ldots)\}$ and $\mathcal{P}(\eta)=span\{E_{01}\oplus 0,E_{10}\oplus 0,(3E_{00}-E_{11})\oplus 0, (E_{00}+E_{11})\oplus (1,1,1,\ldots)\}$, where $E_{ij}$'s denote standard matrix units. Let $\circ_\tau$, and $\circ_\eta$ be the Choi-Effros products on the peripheral Poisson boundaries $\mathcal{P}(\tau)$ and $\mathcal{P}(\eta)$ respectively. Since $(E_{01}\oplus0)\circ_\tau (E_{10}\oplus0)=E_{00}\oplus (\frac{1}{2},\frac{1}{2},\frac{1}{2},\ldots)\neq E_{00}\oplus 0$, and $(E_{01}\oplus0)\circ_\tau (E_{10}\oplus0)=E_{00}\oplus (\frac{1}{4},\frac{1}{4},\frac{1}{4},\ldots)\neq E_{00}\oplus 0$, we realize that neither $\tau$, nor $\eta$ is peripherally automorphic.
\end{example}

The following special case of Corollary \ref{bell} is of particular interest.
Recall that in finite dimensions normalized trace is a stationary state for mixed unitary maps.

\begin{remark}\label{woods}
Let  $U_i, 1\leq i\leq n$ be some unitaries in $\mathbb{M}_d$. Suppose $\alpha , \beta $ are two proper convex combinations of corresponding automorphisms, that is, $$\alpha (X)=\sum_{i=1}^n p_iU_i^*XU_i, ~~\sum_{i=1}^n p_i=1, ~0< p_i < 1, \ \forall i;$$ and
 $$\beta (X)=\sum_{i=1}^n q_iU_i^*XU_i, ~~ \sum_{i=1}^n q_i=1, ~0<q_i<1, \forall i.$$ Then from the  Corollary \ref{bell} above, it is clear that $\mathcal{P}(\alpha )=\mathcal{P}(\beta)$. In fact, they share the same peripheral eigenvalues and the corresponding eigenspaces also coincide, that is, $E_\lambda(\alpha )=E_\lambda(\beta)$ for all $ \lambda\in\mathbb{T}$.
\end{remark}

In Corollary \ref{bell}, typically there is proper containment in  $\mathcal{P}(\tau )\subseteq  \bigcap _{i=1}^n\mathcal{P}(\tau _i).$
For instance, for every unitary channel on $\mathbb{M}_d$, the peripheral Poisson boundary is the whole algebra $\mathbb{M}_d$ (Example 3.1, \cite{BKT1}).
However, for  mixed unitary channels which are not unitary channels, the peripheral Poisson  boundary must be a proper subalgebra of $\mathbb{M}_d$, as the UCP map must act as an automorphism on the peripheral Poisson boundary.

Stationarity is a crucial assumption in most of these results. The following example shows that the picture can be very different, even for mixed unitary maps acting on the infinite dimensional von Neumann algebras.

\begin{example}\label{t,1-t}
Consider the Hilbert space $\mathcal{H}=L^2 [0,1]$. Let $g:[0,1]\to [0,1]$ be a strictly increasing $C^1$- function, such that $g(0)=0$ and $g(1)=1$. Then, by the change of variable formula, $U_g : \mathcal{H}\to \mathcal{H}$, defined by $(U_g f)(t)=f(g(t))\sqrt{g'(t)}$, $\forall t \in [0,1]$, is a unitary operator with $(U_g)^*=U_{g^{-1}}$.

For $v \in L^{\infty}[0,1]$, let $M_v$ be the multiplication operator defined by $(M_v f)(t)=v(t)f(t)$.  We identify $L^\infty[0,1]$ with the von Neumann algebra $\{M_v: v\in L^\infty[0,1]\} \subset \mathscr{B}(L^2[0,1])$ by identifying $v$ with $M_v.$
Observe that $U_g M_v {U_g^*} = M_{v \circ g}$ for any $g$ as above.

Let $g_1 (t)= t^2 \ \mbox{and} \ g_2 (t)= 2t - t^2$, which means that ${g_1}^{-1}(t)=\sqrt{t} \ \mbox{and} \ {g_2}^{-1}(t)= 1- \sqrt{1-t}$. Consider the linear map $\tau: L^{\infty}[0,1] \to L^{\infty}[0,1]$, defined as $\tau(M_v)= \frac{1}{2}(U_{g_1} M_v U^*_{g_1} + U_{g_2} M_v U^*_{g_2})$. Hence, $\tau$ is a mixed unitary map, such that $\tau(M_v)= \frac{1}{2}(M_{v_1} + M_{v_2})$, where $v_1 (t)= v \circ g_1 (t) \ \mbox{and} \ v_2 (t)= v \circ g_2 (t)$. 
For $v_0(t)=t$,  $\tau(M_{v_0})= M_{v_0}$. 
Therefore, $M_{v_0}$ is a fixed point of $\tau$.
$M_{v_0} \geq U_{g_1} M_{v_0} U^*_{g_1}= M_{g_1} \ \mbox{since } \ t\geq t^2$, and $M_{v_0} \leq U_{g_2} M_{v_0} U^*_{g_2}= M_{g_2} \ \mbox{since}, \ t \leq 2t-t^2$. 
Hence, though $M_{v_0}$ is a fixed point of $\tau$, it is neither a fixed point of the mapping $M_v \mapsto U_{g_1}M_v {U^*_{g_1}}$, nor of the mapping $M_v \mapsto U_{g_2}M_v U_{g_2}^*$.
\end{example}

We are interested to know about the Choi-Effros product of $M_t$ with itself in the above example. But before that, we present a much needed result to compute the same.

\begin{lemma}\label{elle}
Let $\mathcal{A}\subseteq\mathscr{B}(\mathcal{H})$ be a von Neumann algebra, and let $\tau:\mathcal{A}\to\mathcal{A}$ be a normal UCP map. 
Let $m$  be a natural number and let $x_1, \ldots , x_m\in \mathcal{A}$, such that $\tau(x_j)=\lambda_jx_j,~ \lambda _j\in \mathbb{T},~\forall 1\leq j\leq m$. Then  $$x_1\circ x_2\cdots \circ x_m=\mbox{s.}\lim_{n\to\infty}\frac{\tau^n(x_1x_2\ldots x_m)}{(\lambda_1\lambda_2\cdots \lambda _m)^n}.$$
\end{lemma}
\begin{proof}
By Theorem 1.1 in \cite{BKT1}, for the triple $(\mathcal{H}, \mathcal{A}, \tau )$, there exists a triple $(\mathcal{K}, \mathcal{B}, \theta )$ called  the minimal dilation of $\tau$, such that $\tau^n(x)=p\theta^n(x)p, ~\forall x\in\mathcal{A}$, where $p$ is a projection from $\mathcal{K}$ onto $\mathcal{H}$. Moreover, by Lemma 2.1 in \cite{BKT1}, we have that for $1\leq j\leq m$, there exist unique $\widehat{x}_j \in\mathcal{B}$, such that $\theta(\widehat{x}_j)=\lambda_j\widehat{x}_j$ and $x_j=p\widehat{x}_jp.$ It is given by
$$\widehat{x}_j= \mbox{s.}\lim _{n\to \infty}\lambda ^{-n}\theta ^n(x_j),$$
and $x_1\circ x_2\circ \cdots \circ x_m= p(\widehat{x}_1\widehat{x}_2\ldots \widehat{x}_m)p.$  Hence, as products of operators  is jointly continuous in SOT for sequences,

$$x_1\circ x_2\circ \cdots \circ x_m= \mbox{s.}\lim _{n\to \infty}\frac{p\theta ^n(x_1x_2\ldots x_m)p}{(\lambda _1\lambda _2\ldots \lambda _m)^n}= \mbox{s.}\lim _{n\to \infty}\frac{\tau ^n(x_1x_2\ldots x_m)}{(\lambda _1\lambda _2\ldots \lambda _m)^n}.$$
\end{proof}

\begin{remark}
It has been observed that the functions $1_{[0,1]}, M_t \in \mathcal{P}(\tau)$. We want to show that $M_t$ is a projection in the peripheral Poisson boundary. Clearly it is self-adjoint. So it suffices to show that $M_t\circ M_t=M_t$. By the formula (\ref{limit formula}) for the extended Choi-Effros product, $M_t \circ  M_t=\mbox{s.}\lim_{n\to\infty}\tau^n(M_{t^2})$, and $M_t\geq \tau^n(M_{t^2}), \ n\in\mathbb{N}$, hence, we have $M_t\geq M_t\circ M_t$. Moreover, $M_{t^2}\leq \tau^n(M_{t^2})\leq \tau^{n+1}(M_{t^2})$, since $\tau(M_{t^2})=M_{t^2}+M_{t^2(1-t)^2}\geq M_{t^2}$. 
Since $\tau^n(\tau(M_{t^2})-M_{t^2})=\tau^n(M_{t^2(1-t)^2})$ for every $n\in\mathbb{Z}_+$, $\mbox{s.}\lim_{n\to\infty}\tau^n(M_{t^2(1-t)^2})=0$.  By Kadison-Schwarz inequality, since for every $n\in\mathbb{N}$, $\tau^n$ is a UCP map, $\tau^n(M_t-M_{t^2})=\tau^n(M_{t(1-t)})\leq (\tau^n(M_{t^2(1-t)^2}))^{\frac{1}{2}}$, therefore $\mbox{s.}\lim_{n\to\infty}\tau^n(M_t-M_{t^2})=0$.
Therefore, $M_t=M_t\circ M_t$. Moreover, noting Lemma \ref{elle}, we deduce that $\mbox{s.}\lim_{n\to\infty}\tau^n(M_t^m)=M_t^{\circ m}=M_t$ for every $m\in\mathbb{N}$.

Note that, for any polynomial $p(t)=p_0+p_1 t+\cdots+p_m t^m$, $\mbox{s.}\lim_{n\to \infty}\tau^n(M_{p(t)})=M_{p_0+(p_1+\cdots+p_m)t}$. For $f\in C[0,1]$, we show that $\mbox{s.}\lim_{n\to\infty}\tau^n(M_{f(t)})=f(0)1_{[0,1]}+(f(1)-f(0))M_t$. By Stone-Weierstrass theorem, there exists a sequence of polynomials $\{A_l(t)=a^{(l)}_0+a^{(l)}_1t+\cdots+a^{(l)}_{k_l}t^{k_l}\}_{l\in\mathbb{N}}$, which converges to $f$ in the sup-norm. For $\epsilon>0$, there exists a $m$, such that $||f-A_m||<\frac{\epsilon}{3}$. 
Let $v\in L^2[0,1]$ such that $||v||_2=1$. Certainly, there exists a $n$ such that $||\{\tau^n(M_{A_m(t)})-M_{a^{(m)}_0+[a^{(m)}_1+\cdots+a^{(m)}_{k_m}]t}\}(v)||<\frac{\epsilon}{3}$.
Therefore, $||\{\tau^n(M_{f(t)})-M_{f(0)+t[f(1)-f(0)]}\}(v)||\leq||\tau^n(M_{f(t)})-\tau^n(M_{A_m(t)})||+||\{\tau^n(M_{A_m(t)})-M_{a^{(m)}_0+[a^{(m)}_1+\cdots+a^{(m)}_{k_m}]t}\}(v)||+||(a^{(m)}_0+[a^{(m)}_1+\cdots+a^{(m)}_{k_m}]t)-(f(0)+t[f(1)-f(0)]) ||<\epsilon$. Since $v\in L^2[0,1]$ and $\epsilon>0$ were chosen arbitrarily, we have the claim.

We further claim that if $\tau(M_f)=\lambda M_f$ for $\lambda\in\mathbb{T}$, where $0\neq f\in C[0,1]$, then $\lambda=1$, and $f(t)=a1+bt$ for some  $a,b\in\mathbb{R}$. The situation  for fixed points, that is, when $\lambda =1$ is clear from the computations made above. Now if $\lambda\neq 1$, since $f(0)=\lambda f(0)$ and $f(1)=\lambda f(1)$, implies $f(0)=f(1)=0$. We shall show that $f=0$. If not, then there exists a $c\in(0,1)$, where $|f(c)|=\sup _{t\in[0,1]}|f(t)|$. Since $\frac{1}{2}(f(c^2)+f(2c-c^2))=\lambda f(c)$, by triangle-inequality, $|f(c)|\leq \frac{1}{2}(|f(c^2)|+|f(2c-c^2)|)$, which implies that $|f(c)|=|f(c^2)|$. Therefore, $c^2$ also a point of supremum for $|f|$. Continuing the process, we have $|f(c)|=|f(c^{2^n})|, ~\forall n\in\mathbb{N}$. By continuity of $f$, $f=0$. 
\end{remark}

Now, take a look at the following example.

\begin{example}\label{freegroup}
Consider the discrete group $\Gamma=\mathscr{F}_2$, that is, the free non-abelian group generated by a and b. For the Hilbert space $\mathcal{H}=l^2 (\Gamma)$, let $U: \Gamma \to \mathscr{B}(\mathcal{H})$  be the left regular representation defined as $U_s (\delta_t)=\delta_{st}, \ \forall s,t \in \Gamma$, for the orthonormal basis $\{\delta_t : t \in \Gamma\}$ of $\mathcal{H}$. Consider the following subsets of $\Gamma$, defined by putting elements in their reduced form: 
\begin{eqnarray*}
G_1&=&\{b^{n_k}a^{m_k}\cdots b^{n_2}a^{m_2}b^{n_1}a^{m_1} ; \ m_i, \ n_i \in \mathbb{Z}, \ m_1 <0\},  \\
G'_1&=&\{a^{m_k}\cdots b^{n_2}a^{m_2}b^{n_1}a^{m_1};  \  m_i, \ n_i \in \mathbb{Z}, \ m_1 <0\},  \\
G_2&=&\{a^{n_l}b^{m_l}\cdots a^{n_2}b^{m_2}a^{n_1}b^{m_1}; \  m_i, \ n_i \in \mathbb{Z}, \ m_1 <0\},  \\
G'_2&=&\{b^{m_l}\cdots a^{n_2}b^{m_2}a^{n_1}b^{m_1}; \  m_i, \ n_i \in \mathbb{Z}, \ m_1 <0\},\\
K&=&  \{ b^{n_k}a^{m_k}\cdots b^{n_1}a^{m_1}, \ a^{m_k}\cdots b^{n_1}a^{m_1}, \ a^{n_k}\cdots a^{n_1}b^{m_1}, \ b^{m_k}\cdots a^{n_1}b^{m_1}; \  m_1 \geq 2,  \ m_i, \ n_i \in \mathbb{Z}\},\\
L&=& \{ a^{m_k}\cdots b^{n_1}a^{m_1}b, \ b^{n_k}\cdots b^{n_1}a^{m_1}b, \ a^{n_k}\cdots a^{n_1}b^{m_1}a, \ b^{m_k}\cdots a^{n_1}b^{m_1}a; \  m_1 \geq 1, \ m_i, \ n_i \in \mathbb{Z}\}.
\end{eqnarray*}  
Now define $f \in l^{\infty}(\Gamma)$ by:
$$f(x)=\begin{cases}
1,& ~\text{if}~ x\in\{1, a,  b\}\cup \ G_1\cup \ G'_1\cup \ G_2\cup \ G'_2\cup \ G_2 a\cup \ G'_2 a\cup  \ G_1 b\cup  \ G'_1 b; \\
2, & ~\text{if}~x\in K;\\
0, & ~\text{if}~ x\in L.
\end{cases}$$
Consider the mixed unitary  map $\tau$ defined as $$\tau(M_h)=\frac{1}{2}(U^*_a M_h U_a + U^*_b M_h U_b) , \ \forall h \in l^{\infty}(\Gamma).$$ Note that $U^*_s M_h U_s (t)= M_{h(st)}, \ \forall s,t \in \Gamma$. Elementary calculations show that $\tau(M_f)=M_f$. Since $f(ab)=0, \mbox{and} \ f(b)=1$, $U^*_a M_f U_a (b)=M_{f(ab)} \neq M_{f(b)}$, which implies that $U^*_a M_f U_a  \neq M_{f}$.
\end{example}

The example below illustrates that for normal UCP maps which are convex combinations of the same UCP maps but with different choices of the $p_j$’s, a peripheral eigenvector for one map need not necessarily remain peripheral eigenvector for another, even though this property does hold in finite dimensions.
\begin{example}\label{unitary mix}
Consider the free group $\Gamma=\mathcal{F}_2$ generated by $a$ and $b$, with the Hilbert space $\mathcal{H}=l^2(\Gamma)$, and the von Neumann algebra $\mathcal{A}=l^\infty(\Gamma)$. Let $\tau$ be a normal UCP map with $f\in\mathcal{F}(\tau)$ defined as in the above example. Let the convex combination of unitary maps $\eta$ acting on $\mathcal{A}$ be defined as $\eta(M_h)=\frac{1}{4}(U_a^*M_hU_a+3U_b^*M_hU_b), ~\forall h\in l^\infty (\Gamma)$. Therefore, 
$\eta (M_f)=M_g$, where $g(x )= \frac{1}{4}(f(ax )+3f(bx)).$ However,
$g(a)=\frac{1}{4}(f(a^2)+3f(ba))=\frac{1}{2}\neq 1=f(a)$. Hence, $M_f$ is not a fixed point of $\eta$.
\end{example}

It is clear from  Example \ref{t,1-t} and Example \ref{freegroup} that a result as in Theorem \ref{finite unitary} does not hold in general in infinite dimensions. It turns out that if the unitaries commute mutually, then we have an affirmation. In fact,  a much stronger result holds, valid for commuting contractions on general Banach spaces.  To derive which, we recall a standard version of Choquet-Deny theorem.

\begin{theorem}\label{choquetdeny} \cite{CS} [Choquet-Deny theorem]  
Let $(S,+)$ be an abelian,  topological (Hausdroff) semigroup, and let $\mu$ be a regular Borel probability measure on $S$. Suppose $h: S \to \mathbb{C}$ is a bounded, Borel measurable, complex-valued function such that it is $\mu-$harmonic ($h(x)=\int_{S} h(x+y)\mu(dy) , \ \forall x \in S$), then, for each $x\in S$, $h(x)=h(x+y)$, for almost every $y$ with respect to $\mu$.
\end{theorem}

In ref. \cite{BP}, Bebe Prunaru has provided a short and  elegant proof for the Choquet-Deny Theorem \ref{choquetdeny}. That was for $\mu$-harmonic functions. Here we extend the result for $\lambda-\mu$-harmonic functions (defined in the following theorem) for $\lambda\in\mathbb{T}$.

\begin{theorem}\label{choquetdenyperipheral}
Let $(S,+)$ be an abelian topological (Hausdroff) semigroup with a regular Borel probability measure $\mu$. Let $h:S\to\mathbb{C}$ be a bounded Borel measurable function such that $h$ is $\lambda-\mu$-harmonic for $ \lambda\in\mathbb{T}$, that is,  $$h(x)=\lambda\int_S h(x+y)\mu(dy), \ \forall x\in S.$$ 
Then, for each $x\in S$,  $h(x)=\lambda h(x+y)$ for almost every $y\in S$ with respect to $\mu$.
\end{theorem}

\begin{proof}
Let $\Phi: l^\infty (S)\to l^\infty(S)$ be a bounded linear map defined as $$\Phi(f)(x)=\int_S f(x+y)\mu(dy), \ \mbox{for} \ f\in l^\infty(S).$$
Clearly, $h=\lambda \Phi(h)$. Let $g(x)=\int_S |h(x)-\lambda h(x+y)|^2\mu (dy).$ By expanding the $|\cdot |^2$ in the integral and using $|\lambda |=1$, we see $g=\Phi(|h|^2)-|h|^2$. 
It turns out that $\Phi(g)\geq g$. This can be seen using Fubini's theorem: 
\begin{eqnarray*} 
\Phi(g)(x)&=&\int_S g(x+z)\mu(dz)\\
&=&\int_S \int_S |h(x+z)-\lambda h(x+z+y)|^2\mu (dy)\mu(dz)\\
&=&\int_S \int_S |h(x+z)-\lambda h(x+z+y)|^2\mu (dz)\mu(dy)\\
&\geq & \int_S |\int_S h(x+z)-\lambda h(x+z+y)\mu (dz)|^2\mu(dy)\\
&=&g(x).\end{eqnarray*}
Therefore, $\Phi^{n+1}(|h|^2)-\Phi^n(|h|^2)\geq g\geq 0$. The sequence $\{\Phi^n(|h|^2)\}_{n\in\mathbb{N}}$ converges uniformly, as $||h||_{\infty} ^2 \geq \Phi^{n+1}(|h|^2)\geq \Phi^n(|h|^2)$. Taking $\lim_{n\to \infty}\Phi^{n+1}(|h|^2)-\Phi^n(|h|^2)\geq g\geq 0$  results in $g=0$. Hence, the claim.
\end{proof}

In the above theorem, it is important to note that if $\mu(0)\neq 0$, and $h\neq 0$ in $l^\infty$, then the value of $\lambda$ is forced to be 1. We illustrate the theorem with the following example.
 
\begin{example}
Consider the set of natural numbers $\mathbb{N}$ which is a semigroup under multiplication. Let the probability measure $\mu$ on $\mathbb{N}$ be supported on $\{2\}$, that is, $\mu(\{2\})=1$ and $\mu(\{n\})=0$ for $n\neq 2.$ 
For $\lambda\in\mathbb{T}$, let $h\in l^\infty(\mathbb{N})$ be as follows:
$$h(n)=\begin{cases}
    1, & ~\text{if}~ n ~ \mbox{is odd},\\
    \lambda^k, &~\text{if}~ n=2^k (2m+1) ~\mbox{for }m\in\mathbb{Z}_+.
\end{cases}.$$
With elementary calculations one can  see that $h$ is $\overline{\lambda}-\mu$-harmonic. By the extended Choquet-Deny Theorem above, for every $n\in\mathbb{N}$, $h(n)=\overline{\lambda}h(2n)$. To take a look at this, let's split $n$ for odd and even values. If $n=2m+1$ for $m\in\mathbb{Z}_+$, $h(2m+1)=1=\overline{\lambda}h(2(2m+1))$, and for $n=2^k(2m+1)$ where $k>0$ and $m\in\mathbb{Z}_+$, $h(2^k(2m+1))=\lambda^k=\overline{\lambda}\lambda^{k+1}=\overline{\lambda}h(2^{k+1}(2m+1))$. 
\end{example}

In the next theorem, we make use of multinomial theorem, which is valid for commuting terms.
Suppose $\tau $ is a linear map such that $\tau =\beta _1+\beta _2+\cdots +\beta _n$, where $\beta _1, \ldots , \beta _n$ are mutually commuting, then, 
\begin{equation} \label{binomial}
\tau^N (x)=\sum_{k_1 + \cdots  + k_n =N }\binom{N}{k_1, k_2, \ldots ,k_n}\beta _1^{k_1}\beta _2^{k_2}\cdots \beta _n^{k_n}(x),   \ \forall N \in \mathbb{N},
\end{equation} where $\binom{N}{k_1, k_2, \ldots ,k_n}=\frac{N!}{k_1 !k_2 ! \cdots k_n !}$.


\begin{theorem}\label{Banach}
Let $\mathcal{X}$ be a Banach space and let $\tau : \mathcal{X}\to\mathcal{X}$ be a contractive linear map.  Suppose
$\tau =p_1\alpha _1+p_2\alpha _2+\cdots +p_n\alpha _n$, where $\alpha_1, \alpha_2 , \ldots , \alpha_n $  are contractive linear maps which mutually commute, that is, $\alpha_s\alpha_r = \alpha_r \alpha_s$ for  $1\leq r,s \leq n$, and $\sum_{j=1}^n p_j =1 ~\mbox{such that}~0<p_j<1, ~\forall j$.  If $ \tau(x)=\lambda x$, for some $\lambda \in \mathbb{T}$,  then $\alpha_j (x)= \lambda x$, $\forall j$. 
\end{theorem}
\begin{proof}
Consider the semigroup $S \subseteq \mathbb{Z}_+\times \mathbb{Z}_+\times \cdots \times \mathbb{Z}_+ ~(n+1-\mbox{times}~)$, such that every $s \in S$ is of the form $(R,r_1,r_2,\ldots,r_n)$, where $R= r_1+\cdots+r_n.$ Defining the probability measure as $\mu(R,r_1,r_2,\ldots,r_n)= \frac{p_1^{r_1}p_2^{r_2}\cdots p_n^{r_n}}{2^{R+1}}\frac{R!}{r_1 !r_2 !\cdots r_n !}$.
       
$$\sum_{(N,k_1,k_2,\ldots,k_n)\in S} \mu{(N,k_1,k_2,\ldots,k_n )}=\sum_{(N,k_1,k_2,\ldots,k_n)\in S}{\binom{N}{k_1, k_2, \ldots ,k_n}p_1^{k_1} p_2^{k_2} \cdots p_n^{k_n}}(\frac{1}{2^{N+1}})$$
$$= \sum_{N\geq 0}\frac{1}{2^{N+1}}(\sum_{k_1+k_2+\cdots+k_n= N}{\binom{N}{k_1, k_2, \ldots ,k_n}p_1^{k_1} p_2^{k_2}\cdots p_n^{k_n}})=\sum_{N \geq 0}\frac{1}{2^{N+1}}=1.$$

Let $\rho $ be a bounded linear functional on $\mathcal{X}.$ Define a bounded Borel-measurable function $h: S \to \mathbb{C}$ as  
$$h(N,k_1,k_2,\ldots,k_n)=\rho(\frac{\alpha^{k_1}_1 \alpha^{k_2}_2 \cdots\alpha_n^{k_n}(x)}{\lambda^{N}}),~~ \forall (N,k_1,k_2,\ldots,k_n) \in S.$$
It is easy to verify that $h$ is $\mu$-harmonic. This is because, (using $\tau^M(x) =\lambda^M x,\forall M\in \mathbb{N} $) for any $(R,r_1,r_2,\ldots,r_n)\in S$, 
\begin{eqnarray*}
& &\sum_{(N,k_1,k_2,\ldots,k_n)\in S} h(R+N,r_1+k_1,r_2+k_2,\ldots,r_n+k_n) \mu(N,k_1,k_2,\ldots,k_n)\\&= &\sum_{(N,k_1,k_2,\ldots,k_n) \in S}\!\!\!\!\!\!\!\rho(\frac{\alpha^{k_1}_1 \alpha^{k_2}_2 \cdots \alpha^{k_n}_n (\alpha^{r_1}_1 \alpha^{r_2}_2 \cdots \alpha^{r_n}_n(x))}{\lambda^{N+R}}){\binom{N}{k_1, k_2, \ldots ,k_n}p_1^{k_1} p_2^{k_2} \cdots p_n^{k_n}}(\frac{1}{2^{N+1}})\\ &=& \sum_{N \geq 0} \rho(\frac{\tau^{N}(\alpha^{r_1}_1 \alpha^{r_2}_2 \cdots \alpha^{r_n}_n(x))}{\lambda^{N+R}})(\frac{1}{2^{N+1}})\\&=& \sum_{N \geq 0} \rho(\frac{\alpha^{r_1}_1 \alpha^{r_2}_2 \cdots \alpha^{r_n}_n(x)}{\lambda^{R}})(\frac{1}{2^{N+1}}) \\&=& \rho(\frac{\alpha^{r_1}_1 \alpha^{r_2}_2 \cdots \alpha^{r_n}_n(x)}{\lambda^{R}})\\&=&h(R,r_1,r_2,\ldots,r_n).
\end{eqnarray*} 
Then, by Choquet-Deny theorem (Theorem \ref{choquetdeny}), we get that $$h(0,0,0,\ldots,0)=h(1,1,0,\ldots,0)=h(1,0,1,\ldots,0)=\cdots=h(1,0,0,\ldots,1),$$ which implies that $\rho (\alpha_1(x))=\rho (\alpha_2(x))=\cdots=\rho (\alpha_n(x))=\lambda \rho (x)$.  This completes the proof, since $\rho$ is an arbitrary bounded linear functional.
\end{proof}

The following special case is what we are currently interested in.

\begin{corollary}\label{main}
Let $\mathcal{A}$ be a $C^*$-algebra and $\alpha_1 , \alpha_2 , \ldots , \alpha_n$ be contractive (unital included) completely positive maps on $\mathcal{A}$, such that they mutually commute, that is, $\alpha_s \alpha_r = \alpha_r \alpha_s$ for every $1\leq r,s\leq n$. Define the contractive (unital) completely positive map $\tau: \mathcal{A}\to \mathcal{A}$ as $\tau(y)=p_1 \alpha_1 (y) + p_2 \alpha_2 (y) + \cdots + p_n \alpha_n (y)$, where $\sum_{j=1}^n p_j =1$ and  $0<p_j<1, ~~\forall j.$ If $\tau(x)=\lambda x$, for some $\lambda \in \mathbb{T}$, then $\alpha_j (x)= \lambda x$, for every $j$. 
\end{corollary}

To illustrate this result, we give a simple example  using the Grigorchuk group.

\begin{example}\cite{RGIP}\label{Grigorchuk}        
Consider a binary rooted tree $T_2$, with vertices as the set $\Sigma ^*$ of finite strings in the alphabet $\Sigma =\{0,1\}$, including the empty string $\emptyset $, which roots $T_2$. For $x\in \Sigma  ^*$, the string $x0$ is called the left child of $x$, and the string $x1$ is called the right child of $x$.  The group of all automorphisms $Aut(T_2)$, is a group of all length-preserving permutations $\sigma$ on $\Sigma ^*$, that also respects the initial segment relation: for every string $x$ with initial segment say $y$, then  $\sigma(y)$  will be the initial segment of the string $\sigma(x)$.

The Grigorchuk group $G$ is a subgroup of $Aut(T_2)$, generated by four elements $a ,b ,c ,d$ defined as follows: 
$$a(\emptyset )=b(\emptyset)=c(\emptyset)=d(\emptyset)=\emptyset ;$$
$$a(0)=1,   a(1)=0, ~~ b(0)=c(0)=d(0)=0,~  b(1)=c(1)=d(1)=1.$$
       
For every string $x \in \Sigma ^*$,
\begin{eqnarray*}
a(0x)=1x, \ & a(1x)=0x, \\b(0x)=0a(x), \ & b(1x)=1c(x),\\c(0x)=0a(x), \ & c(1x)=1d(x),
\\d(0x)=0x, \ & d(1x)=1b(x).
\end{eqnarray*} 
The Grigorchuk group $G$ is infinite, where $b,c,d$ mutually commute and $a^2=b^2=c^2=d^2=e$. In fact, $bc=d, ~ cd=b, ~ bd=c$, hence, the group $G=\langle a,b,c \rangle=\langle a,c,d\rangle=\langle a,b,d\rangle$.
Consider the von Neumann algebra $\mathcal{A}:=l^{\infty}(\Sigma  ^*)$. Define the UCP maps $\tau_b,$ and $\tau_c$  as follows:
$$\tau_b(f)(x)=f(b(x)), ~~~~\tau_c(f)(x)=f(c(x)), ~~~~\forall  f\in \mathcal{A}, ~~x \in \Sigma ^*,$$
Observe that $\tau _b$ and $\tau _c$ commute. Take $\tau=\frac{1}{2}(\tau_b + \tau_c).$ 

Now we analyze the peripheral Poisson boundaries of these UCP maps. Let $0\neq f,g\in l^\infty(\Sigma^*)$, such that $\tau_b(f)=\lambda_1 f$ and $\tau_c(g)=\lambda_2 g$ for $\lambda_1,\lambda_2\in\mathbb{T}$, which means $f(b(x))=\lambda_1 f(x)$ and $g(c(x))=\lambda_2 g(x)$ for every string $x\in \Sigma^*$. Since $b^2=c^2=e$, the peripheral eigenvalues of $\tau_b$ and $\tau_c$ are $1$ and $-1$. We have the following: \begin{eqnarray*}
f(\emptyset)=\lambda_1 f(\emptyset)&g(\emptyset)=\lambda_2 g(\emptyset)\\ 
    f(0)=\lambda_1 f(0)&g(0)=\lambda_2 g(0)\\ 
    f(1)=\lambda_1 f(1)&g(1)=\lambda_2 g(1)\\ f(10)=\lambda_1 f(10)&g(10)=\lambda_2 g(10)\\ f(11)=\lambda_1 f(11)&g(11)=\lambda_2 g(11)\\ 
    f(01z)=\lambda_1 f(00z), & ~~~g(01z)=\lambda_2 g(00z),\\
     f(00z)=\lambda_1 f(01z), & ~~~g(00z)=\lambda_2 g(01z),\\
      f(101z)=\lambda_1 f(100z), & ~~~g(10z)=\lambda_2 g(10z),\\
       f(100z)=\lambda_1 f(101z), & \\
        f(110z)=\lambda_1 f(110z), & g(110)=\lambda_2 g(110)\\
        &~~~g(1100z)=\lambda_2 g(1101z),\\
          & ~~~g(1101z)=\lambda_2 g(1100z),\\
          f(111b(z))=\lambda_1 f(111z), & ~~~g(111c(z))=\lambda_2 g(111z),
\end{eqnarray*} for every string $z\in\Sigma ^*$. For any string $x\in\Sigma ^*$, let $\delta _{x}$ denote the characteristic function of the set $\{x\}$, considered as an element of $\mathcal{A}.$
 Below $(x)^ n$ will denote $x$ repeated $n$ times, therefore $(x)^ 0=\emptyset$ and $(x)^ 0 y=y$ for every $x,y\in\Sigma^*$; and $\overline{\mathcal{Y}}^{w^*}$ denotes the $\mbox{weak}^*$-closure of the subspace $\mathcal{Y}\subseteq l^\infty (\Sigma ^*)$. 
Then the eigenspaces of $\tau_b$ and $\tau_c$ can be described as follows.

\begin{align*}
E_{-1}(\tau_b) &= \overline{\text{span}}^{w^*}\{\delta_{(111)^n 01z}-\delta_{(111)^n 00z}, \delta_{(111)^n 101z}-\delta_{(111)^n 100z} : z\in\Sigma^*, n\in\mathbb{Z}_+\},\\
E_1(\tau_b) &= \overline{\text{span}}^{w^*}\{\delta_{(1)^n}, \delta_{(1)^n 0}, \delta_{(111)^n 110z},\\
&\qquad \delta_{(111)^n 01z}+\delta_{(111)^n 00z}, \delta_{(111)^n 101z}+\delta_{(111)^n 100z} : z\in\Sigma^*, n\in\mathbb{Z}_+\};\\
E_{-1}(\tau_c) &= \overline{\text{span}}^{w^*}\{\delta_{(111)^n 01z}-\delta_{(111)^n 00z}, \delta_{(111)^n 1100z}-\delta_{(111)^n 1101z} : z\in\Sigma^*, n\in\mathbb{Z}_+\},\\
E_1(\tau_c) &= \overline{\text{span}}^{w^*}\{\delta_{(1)^n}, \delta_{(1)^n 0}, \delta_{(111)^n 01z}+\delta_{(111)^n 00z},\\
&\qquad \delta_{(111)^n 10z}, \delta_{(111)^n 1100z}+\delta_{(111)^n 1101z} : z\in\Sigma^*, n\in\mathbb{Z}_+\}.
\end{align*}

Hence, the peripheral Poisson boundaries are $\mathcal{P}(\tau_b)=\mathcal{P}(\tau_c)=l^\infty(\Sigma  ^*)$, where the Choi-Effros product coincides with regular multiplication. Next, we compute the peripheral Poisson boundary of $\tau$. If $\tau(h)=\lambda h$ for $h\neq0,~\lambda\in \mathbb{T}$, then by Corollary \ref{main}, $\lambda\in\{-1,1\}$, and $\tau_b(h)=\lambda h$ and $\tau_c(h)=\lambda h$. Using $\frac{1}{2}(h(b(x))+h(c(x)))=\lambda h(x)$ for every string $x\in\Sigma ^*$, we deduce the peripheral eigenspaces as:

\begin{align*}
E_{-1}(\tau)&=\overline{span}^{w^*}\{\delta_{(111)^ n00z}-\delta_{(111)^ n01z}: z\in\Sigma ^*,n\in\mathbb{Z}_+\},\\  
E_1(\tau)&=\overline{span}^{w^*}\{\delta_{(1)^ n},\delta_{(1)^ n0},\delta_{(111)^ n00z}+\delta_{(111)^ n01z},\\
& \qquad \delta_{(111)^ n100z}+\delta_{(111)^ n101z}, \delta_{(111)^ n1100z}+\delta_{(111)^ n1101z}: z\in\Sigma ^*,n\in\mathbb{Z}_+\}.
\end{align*}

The peripheral Poisson boundary is 

\begin{align*}
\mathcal{P}(\tau)&=\overline{span}^{w^*}\{\delta_{(1)^ n},\delta_{(1)^ n0},\delta_{(111)^ n00z},\delta_{(111)^ n01z},\delta_{(111)^ n100z}+\delta_{(111)^ n101z},\\ &\qquad \delta_{(111)^ n1100z}+\delta_{(111)^ n1101z}: z\in\Sigma ^*,n\in\mathbb{Z}_+\},\end{align*} with the Choi-Effros product coinciding with the regular multiplication.    
\end{example}

In  Corollary \ref{main},  we have that if $x\in E_\lambda (\tau)$, for some $\lambda \in \mathbb{T}$, then $x\in E_{\lambda }(\alpha _j)$, for every $j$. Therefore, it is clear that $\mathcal{P}(\tau )\subseteq \bigcap _j\mathcal{P}(\alpha _j)$ as subspaces of $\mathcal{A}.$
Hence, we wonder if $\mathcal{P}(\tau )= \bigcap _j\mathcal{P}(\alpha _j)$. In general, this is not true. We illustrate this with the following example.

\begin{example}
Consider the UCP maps $\alpha,\beta:\mathbb{M}_3\to\mathbb{M}_3$ defined as follows: 
$$\alpha(\begin{bmatrix}
x_{11}&x_{12}&x_{13}\\x_{21}&x_{22}&x_{23}\\x_{31}&x_{32}&x_{33}
\end{bmatrix})=\begin{bmatrix}
x_{11}&(\frac{\sqrt{3}+i}{2})x_{12}&(\frac{1+i\sqrt{3}}{2})x_{13}\\(\frac{\sqrt{3}-i}{2})x_{21}&x_{22}&(\frac{\sqrt{3}+i}{2})x_{23}\\(\frac{1-i\sqrt{3}}{2})x_{31}&(\frac{\sqrt{3}-i}{2})x_{32}&x_{33}
\end{bmatrix},$$ and 
$$\beta(\begin{bmatrix}
x_{11}&x_{12}&x_{13}\\x_{21}&x_{22}&x_{23}\\x_{31}&x_{32}&x_{33}
\end{bmatrix})=\begin{bmatrix}
x_{11}&ix_{12}&-x_{13}\\-ix_{21}&x_{22}&ix_{23}\\-x_{31}&-ix_{32}&x_{33}
\end{bmatrix}.$$ 
    
Clearly, $\alpha$ and $\beta$ mutually commute. Let $\tau:\mathbb{M}_3\to\mathbb{M}_3$ be a UCP map defined as $\tau=\frac{\alpha+\beta}{2}$, i.e.

$$\tau(\begin{bmatrix}
x_{11}&x_{12}&x_{13}\\x_{21}&x_{22}&x_{23}\\x_{31}&x_{32}&x_{33}\end{bmatrix})=\begin{bmatrix}
x_{11}&(\frac{\sqrt{3}+3i}{4})x_{12}&(\frac{-1+i\sqrt{3}}{4})x_{13}\\(\frac{\sqrt{3}-3i}{4})x_{21}&x_{22}&(\frac{\sqrt{3}+3i}{4})x_{23}\\(\frac{-1-i\sqrt{3}}{4})x_{31}&(\frac{\sqrt{3}-3i}{4})x_{32}&x_{33}
\end{bmatrix}.$$

Clearly, the peripheral Poisson boundaries for $\tau,\alpha,\beta$ are $\mathcal{P}(\tau)=\mathbb{C}\oplus\mathbb{C}\oplus\mathbb{C}$, and $\mathcal{P}(\alpha)=\mathbb{M}_3=\mathcal{P}(\beta)$.
\end{example}

We are interested to know for $\mathbb{M}_d$, if we drop the condition of $\alpha _j$'s commuting mutually, then, would the conclusion of Corollary \ref{main} still hold?
The answer is no. We illustrate this with the following example.

\begin{example}
Consider the unital $C^*$-algebra  $\mathcal{A}=\mathbb{C}^3$, and let $\tau_1,\tau_2 : \mathcal{A}\to \mathcal{A}$ be UCP maps defined as follows:
$\tau_1 (x_1 , x_2 , x_3) = (x_1 , x_2 , x_2)$ and $\tau_2 (x_1 , x_2 , x_3) = (x_1 , x_2 , x_1)$. Clearly, $\tau_1\tau_2\neq \tau_2\tau_1$. Define the UCP map $\tau=\frac{1}{2}(\tau_1 + \tau_2)$, that is, $\tau (x_1 , x_2 , x_3)= (x_1 , x_2 , \frac{1}{2}(x_1 + x_2))$.  Note that  $\tau (1 , -1 , 0) = (1 , -1 , 0) $, $\tau_1 ( 1,-1,0 ) = (1 , -1 , -1)$,  and $\tau_2 ( 1,-1,0 ) = ( 1, -1 , 1)$. Therefore, there exists a peripheral eigenvector of $\tau$, which is neither a peripheral eigenvector of $\tau_1$, nor of $\tau_2$. 
\end{example}

Now we ask a more subtle question. Consider a normal UCP map $\tau $ acting on a von Neumann algebra $\mathcal{A}$, which is a convex combination of mutually commuting
normal UCP maps $\alpha _1, \alpha _2, \ldots , \alpha _n$. We have seen that
$$\mathcal{P}(\tau ) \subseteq \bigcap _{j=1}^n\mathcal{P}(\alpha _j),$$
as subspaces of $\mathcal{A}$. The peripheral Poisson boundaries $\mathcal{P}(\tau )$ and $\mathcal{P}(\alpha _j)\mbox{'s}, ~\forall j$, are $C^*$-algebras under their respective extended Choi-Effros product. 
So does this mean that we are having an inclusion of $C^*$-algebras? Surprisingly the answer is no. In general, the extended Choi-Effros product between  two peripheral eigenvectors with respect to $\tau$, need not coincide with that of the constituent $\alpha _j$'s. Here is an
example to illustrate this.

\begin{example}\label{subtle}
Consider the von Neumann algebra $\mathcal{A}=\mathbb{M}_2\oplus l^{\infty}(\mathbb{Z}_+)$, and let  $\tau_1, \tau_2 : \mathcal{A}\to \mathcal{A}$ be normal UCP maps defined as $$\tau_1(\begin{bmatrix} x_{00}& x_{01}\\ x_{10}& x_{11}\end{bmatrix}\oplus (y_0 , y_1,\ldots , y_n , \ldots ))=\begin{bmatrix} x_{00}& x_{01}\\ x_{10}& x_{11}\end{bmatrix}\oplus (x_{00}, x_{11}, y_0 , y_1, \ldots , y_n , \ldots )$$ and $\tau _2$ is the identity map.

Clearly $\tau _1, \tau _2$ commute. Let $\tau:\mathcal{A}\to\mathcal{A}$ be defined as  $\tau=\frac{1}{2}(\tau_1 + \tau_2)$. Since $\tau, \tau_1, \tau_2$ are normal UCP maps, they have their respective peripheral Poisson boundaries. Only $1$ happens to be the peripheral eigenvalue for each of these maps. The peripheral Poisson boundaries are as follows: 
$$\mathcal{P}(\tau)=~\mbox{span}~\{E_{0,0} \oplus (1,0,1,0,\ldots), E_{1,1} \oplus(0,1,0,1,\ldots), E_{1,0}\oplus 0, E_{0,1} \oplus 0 \}= \mathcal{P}(\tau_1)$$ and $\mathcal{P}(\tau_2)= \mathcal{A}$. 
We denote $\circ_\tau, \circ_{\tau_1}, \circ_{\tau_2}$ as the Choi-Effros products of $\tau,\tau_1,\tau_2$ respectively.
Using Formula \ref{limit formula}, it is not hard to see that,
$$(E_{01}\oplus 0)\circ_\tau (E_{10}\oplus 0)=(E_{01}\oplus 0)\circ_{\tau_1} (E_{10}\oplus 0)=E_{00}\oplus(1,0,1,0,\ldots),$$ but $(E_{01}\oplus0)\circ_{\tau_2} (E_{10}\oplus0)=E_{00}\oplus 0$, where $E_{ij}$'s denote matrix units of $\mathbb{M}_2.$
It follows that $\mathcal{P}(\tau )$ and $\mathcal{P}(\tau_1)$ are isomorphic to $\mathbb{M}_2$ and $\mathcal{P}(\tau _2)$ is the same as $\mathcal{A}$.
\end{example}

\section{Continuous time-independence of the peripheral Poisson boundary}
It is shown by Bhat, Kar and Talwar (Theorem 4.3, \cite{BKT1}) that if $\tau $ is  a normal, unital completely positive map  acting on a von Neumann algebra $\mathcal{A}\subseteq \mathscr{B}(\mathcal{H})$, then $\mathcal{P}(\tau)=\mathcal{P}(\tau^k), \ k\in\mathbb{N}$.
Here we extend this result to the  continuous-time case for finite dimensional $C^*$-algebras, using the Jacobs-de Leeuw-Glicksberg theory mentioned before.
 

\begin{theorem}\label{ppbQ}
Let $\mathcal{A}$ be a finite dimensional $C^*$-algebra and let $\tau=\{\tau_t\}_{t\in\mathbb{R}_+}: \mathcal{A} \to \mathcal{A}$ be a quantum dynamical semigroup. 
Then, $\mathcal{P}(\tau)=\mathcal{P}(\tau_t), \forall t>0$, i.e. the peripheral Poisson boundary $\mathcal{P}(\tau _t)$ is independent of $t$ for $t>0,$ and is equal to the peripheral Poisson boundary $\mathcal{P}(\tau )$ of the semigroup.
\end{theorem}
\begin{proof}
Due to finite dimensionality of $\mathcal{A}$, it is reflexive and  $\tau=\{\tau_t\}_{t\in\mathbb{R}_+}$ is almost periodic.
Using the JdLG decomposition theorem, we have
$\mathcal{A}=\mathcal{A}_r(\tau)\oplus\mathcal{A}_{fl}(\tau)$ , where $\mathcal{A}_r(\tau)=\overline{\mbox{span}}\{x\in\mathcal{A}\ : \ y\in \overline{\{\tau_t x\}}_{t\in\mathbb{R}_+} \implies x \in \overline{\{\tau_t y\}}_{t\in\mathbb{R}_+}\}$ and $\mathcal{A}_{fl}(\tau)=\{x\in \mathcal{A} \ : \ \lim_{t\to\infty}\tau_t x=0\}$.  Note that for every character $\alpha: \mathbb{R}_+ \to \mathbb{T}$ corresponding to some $x\in\mathcal{A}_r(\tau)$, there exists an $a\in \mathbb{R},$ such that $\alpha(t)=\exp(iat), \ \forall t$. This is because $\alpha$ is continuous on $\mathbb{R}_+$, as $\tau$ is a ${C}_0$ semigroup. Hence, $\mathcal{A}_r(\tau)=\mathcal{P}(\tau)$.

We show that $P(\tau_t)\subseteq \mathcal{A}_r(\tau), \ \forall t>0$.
Let $x\in E_\lambda(\tau_{t_0})$, then $x \in E_{\lambda^n}(\tau_{n{t_0}}), \ \forall n\in \mathbb{N}$. 
Let $y \in \overline{\{\tau_t x\}}_{t\in\mathbb{R}_+}$.  Now there are two possibilities. Case (i):  $y=\tau_s x$ for some $s\geq 0$. 
Choose a natural number $N>\frac{s}{t_0}.$ Then taking $t_1=Nt_0-s$, we see that $\tau_{t_1}(y)=\tau _{Nt_0}(x)= \lambda^N x$. By the same argument, $\lambda^N x,\lambda^{N+1} x,\lambda^{N+2} x,\ldots \in \{\tau_t y\}_{t\in\mathbb{R}_+}$. We know that there exists a subsequence of $\{\lambda^{N+k}\}_{k\in\mathbb{N}}$ converging to 1. Hence $x\in \overline{\{\tau_t y\}}_{t\in\mathbb{R}_+}$.
Case (ii): $y \neq \tau_s (x), \ \forall s$, but $y\in \overline{\{\tau _t x\}}_{t\in \mathbb{R}_+}$. This implies that there exists a strictly increasing sequence $t_n$, such that $\tau_{t_n}x \to y$ as $n\to \infty$. 
Note that there exists an increasing sub-sequence of natural numbers $\{N_n\}_{n\in \mathbb{N}}$, such that $ N_n t_0 \geq t_n$ for every $n$, and $\lambda^{N_{n}}\to 1$ as $n$ goes to infinity. Then, \begin{equation*}
||\lambda^{N_{n}}x-\tau_{{t_0}N_{n}-t_{n}}(y)||=||\tau_{{t_0} N_{n}}(x)-\tau_{{t_0}N_{n}-t_{n}}(y)||=||\tau_{{t_0}N_{n}-t_{n}}(\tau_{t_{n}}(x)-y)||\leq M_\tau ||\tau_{t_{n}}(x)-y||.
\end{equation*} 
Hence, $x\in  \overline{\{\tau_t y\}}_{t\in\mathbb{R}_+}$. 
Thus, $P(\tau_t)\subseteq \mathcal{A}_r(\tau), \ t>0$. The converse statement $\mathcal{A}_r(\tau)\subseteq P(\tau_t)$ is obvious. 
\end{proof}

We wish to have the same result for adjoints of quantum Makov semigroups.


\begin{corollary}\label{ppbT}
Let $\mathcal{A}$ be a finite dimensional $C^*$-algebra, and $\tau=\{\tau_t\}_{t\in\mathbb{R}_+}:\mathcal{A}\to\mathcal{A}$ be a completely positive trace preserving $\mathcal{
C}_0$ semigroup. Then, $\mathcal{P}(\tau)=\mathcal{P}(\tau_t)$ for every $t>0$.
\end{corollary}
\begin{proof} 
It is to be noted that in the proof of Theorem \ref{ppbQ}, we do not quite need contractivity of the semigroup. As long as the semigroup is a bounded $C_0$-semigroup, we can appeal to Lemma  \ref{Reflex} and the proof would be complete.

By assumption, there exist $d_1,d_2,...,d_m\in\mathbb{N}$, such that $\mathcal{A}=\mathbb{M}_{d_1}\oplus\mathbb{M}_{d_2}\oplus...\oplus\mathbb{M}_{d_m}$. Let $d=d_1+d_2+...+d_m$.
Now complete positivity implies that the norm of $\tau _t$ attained at the identity. As $\text{trace}(\tau _t(1))=\text{trace}(1)=d,$  we can see that $\|\tau _t\|\leq \|\tau _t(1)\|\leq d.$


\end{proof}

Previously in Theorem \ref{pcm}, for finite dimensional $C^*$-algebras, we established a relation between the peripheral Poisson boundary of UCP map, which is a composition of mutually commuting UCP maps, and the peripheral Poisson boundaries of the constituent ones. We show that this relationship carries over to the continuous-time regime as well, that is, for quantum Markov semigroups.

\begin{theorem}\label{continuouscomposition} Consider quantum Markov semigroups $\tau=\{\tau_t\}_{t\in\mathbb{R}_+}$ and 
$ \tau^{(j)}=\{\tau^{(j)}_t\}_{t\in\mathbb{R}_+},$ $ 1\leq j\leq m,$ acting on a finite dimensional $C^*$-algebra $\mathcal{A}$,
with generators $\mathcal{L},\mathcal{L}_j, 1\leq j\leq m$ respectively, such that they all mutually commute, and $\mathcal{L}=\mathcal{L}_1+\mathcal{L}_2+\cdots+\mathcal{L}_m$. Then, $ \mathcal{P}(\tau)=\bigcap_{i=1}^m\mathcal{P}(\tau^{(i)})$.     
\end{theorem}
\begin{proof}  
The quantum Markov semigroups $\tau=\{\tau_t\}_{t\in\mathbb{R}_+}, \tau^{(j)}=\{\tau^{(j)}_t\}_{t\in\mathbb{R}_+}, 1\leq j\leq m,$ are of the form, 
$\tau_t=\exp(t\mathcal{L})$ and $\tau^{(j)}_t=\exp(t\mathcal{L}_j), \ \forall 1\leq j\leq m, \ t\geq 0$. By assumption, $\tau_t=\tau^{(1)}_t\tau^{(2)}_t\cdots\tau^{(m)}_t$ and  by Theorem \ref{ppbQ}, $\mathcal{P}(\tau)=\mathcal{P}(\tau_t),  \  \mathcal{P}(\tau^{(i)})=\mathcal{P}(\tau^{(i)}_t), \ \forall i$ and $t>0$. Using this we deduce that $\mathcal{P}(\tau)=\mathcal{P}(\tau_t)=\bigcap_{j=1}^m\mathcal{P}(\tau^{(j)}_t) = \bigcap_{j=1}^m\mathcal{P}(\tau^{(j)})$.
\end{proof}

We naturally wonder if the time independence follows, when we take into consideration an infinite dimensional von Neumann algebra. This is in fact not true in general, even for abelian von Neumann algebras and semigroups of automorphisms.
 
\begin{example}\label{continuouscounter}
Consider the von Neumann algebra $\mathcal{A}=L^\infty(\mathbb{T})\subseteq\mathscr{B}(L^2(\mathbb{T}))$, and the quantum Markov semigroup $\tau=\{\tau_t\}_{t\in \mathbb{R}_+}$ defined as $\tau_t(f)(z)=f(\exp{(2\pi it)}z), \ \forall t\in \mathbb{R}_+$. We show that, for every $a\in\mathbb{R}_+$
$$\mathcal{P}(\tau _a)=\begin{cases}
C(\mathbb{T})&\text{if} ~a\notin \mathbb{Q};\\
L^{\infty}(\mathbb{T})&\text{if}~a\in \mathbb{Q}.
\end{cases} $$


 First, we consider the irrational  case. Let $a\notin \mathbb{Q}$ and let $f\in\mathcal{A},$ be such that $\tau_a(f)=\mu f$ for some $\mu\in\mathbb{T}$, which means that $f(\exp{(2\pi ia)}z)=\mu f(z)$ for all $z$.  Note that $L^\infty(\mathbb{T})\subseteq L^2(\mathbb{T})$. Let $f(z)=\sum_{n\in\mathbb{Z}}c_nz^n$, be the Fourier series expansion of $f$, where $c_n=\int_\mathbb{T}f(z)z^{-n}dz$ (the integral is with respect to the normalized Haar measure).  Since, $f(\exp{(2\pi ia)}z)=\sum_{n\in\mathbb{Z}}c_n\exp{(2\pi ina)} z^n=\sum_{n\in\mathbb{Z}}c_n \mu z^n=\mu f(z)$, therefore $c_n(\exp{(2\pi ina)}-\mu)=0,
~\forall n\in \mathbb{Z}$. With $a$ being irrational, $\exp{(2\pi ina)}$ does not take the same value for two different $n$'s. This implies there exists an $m\in \mathbb{Z}$, where $\mu=\exp{(2\pi ima)}$ and $c_n=0, \ \forall n\neq m$. Hence, $f(z)=c_mz^m, \forall z\in\mathbb{T}$. Therefore, by Stone-Weierstrass theorem, $\mathcal{P}(\tau_a)=\overline{span}\{z^m:m\in\mathbb{Z}\}=C[0,1]$.

Now if $a\in \mathbb{Q}$ and $a>0$, then $a=\frac{l}{m}$ for some $l,m\in \mathbb{N}$ with g.c.d $(l,m)=1$. In this case, eigenvalues of $\tau _a$ are different $m^{\text{th}}$ roots of unity. Moreover $L^{\infty}(\mathbb{T})$ is seen to be a finite linear span of corresponding eigenspaces. This should also be clear from the discrete situation as well, as $\tau _{am}$ is the identity map.
\end{example}

\section{Concrete structure of the peripheral Poisson boundary  in continuous time}

Based on earlier literature ( See \cite{Lind} and \cite{WP}) a structure theorem for cycles  of trace preserving Schwarz maps is presented in the popular Lecture Notes of Wolf (Theorem 6.16, \cite{WMM}), and the result is well-known.  Building on this, a concrete description of the   peripheral Poisson boundary of discrete semigroups on $\mathbb{M}_d$ was given by D. Amato, P. Facchi, and A. Konderak in  Theorem 3 of \cite{DA}. In this Section we recall this structure theorem and then extend it to continuous time  quantum Markov semigroups. The structure becomes simpler now, as the permutations between different summands is no longer present.


The matrix algebra $\mathbb{M}_d$ can be endowed with a Hilbert space structure by means of the Hilbert-Schmidt scalar product $\langle X,Y\rangle_{HS}= \text{tr}(X^*Y), \forall X,Y\in \mathbb{M}_d $. 
Given a completely positive map $\tau:\mathbb{M}_d \to \mathbb{M}_d$, it's Hilbert space adjoint $\tau^*$ is a completely positive map (Proposition 2.18 in \cite{Watrous}), and obeys $\langle X,\tau(Y)\rangle=\langle\tau^*(X),Y\rangle, \ \forall X,Y\in \mathbb{M}_d$. 
In particular, $\tau:\mathbb{M}_d\to\mathbb{M}_d$ is a unital completely positive map
if and only if $\tau^*$ is a completely positive trace preserving (CPTP) map (Theorem 2.26 in \cite{Watrous}). 
First we look at a special case.

\begin{definition}\label{stationary}
A unital completely positive map $\tau:\mathbb{M}_d \to \mathbb{M}_d$ is said to be {\em stationary}, if it admits a faithful invariant  state: There exists a faithful state $\phi $ of $\mathbb{M}_d$, such that
$$\varphi (\tau (X))= \varphi (X), ~~\forall X\in \mathbb{M}_d.$$
\end{definition}
It is to be noted that $\tau $ is stationary iff
$\tau^*$ has a strictly positive fixed point, that is, there exists $\rho >0$ with $\tau^*(\rho)=\rho $, as any state $\varphi $ on $\mathbb{M}_d$ is given by $\varphi(X)=\text{trace}(\rho X)$ for some density matrix $\rho$. For this reason such maps were called as {\em faithful} maps in \cite{DA}. We do not use that terminology as faithful states has a different meaning.

It is known that if $\tau $ is stationary then it is peripherally automorphic. An explicit description of the peripheral Poisson boundary of such maps was proved in 
\cite{DA}, and we present it here.

\begin{theorem}\cite{DA}
Let $\tau:\mathbb{M}_d\to\mathbb{M}_d$ be a stationary UCP map. Then, there exists a decomposition of the Hilbert space $\mathcal{H}=\mathbb{C}^d$, of the form $\mathcal{H}=\bigoplus_{k=1}^m\mathcal{H}_{k,1}\otimes\mathcal{H}_{k,2}$, where $\mathcal{H}_{k,1}$, $\mathcal{H}_{k,2}$, are Hilbert spaces, such that $\mathcal{P}(\tau)=\bigoplus_{k=1}^M\mathscr{B}(\mathcal{H}_{k,1})\otimes 1_{k,2}$. Therefore, every element of $\mathcal{P}(\tau)$ is of the form $X=\bigoplus_{k=1}^m X_k \otimes 1_{k,2}$, and the action of $\tau$ on $\mathcal{P}(\tau)$ is given by $$\tau(X)=\bigoplus_{k=1}^m (U_{\sigma(k)}^*X_{\sigma(k)}U_{\sigma(k)} \otimes 1_{k,2}),$$ where all the $U_k\mbox{'s, for } 1\leq k\leq m$, are unitaries, and $\sigma$ is a permutation on $\{1,2,\ldots, m\}$ such that the $\dim (\mathcal{H}_{k,1})=\dim (\mathcal{H}_{\sigma(k),1})$. 
\end{theorem}

In other words, the peripheral Poisson boundary is realized as a direct sum of matrix algebras, where $\tau $ acts as permutations between summands having equal dimensions, followed by unitary conjugations.
Now we extend this result to quantum Markov semigroups. The stationarity is defined analogous to the discrete case, by requiring a faithful invariant state. The additional feature we see here is that the CP semigroup $\tau _t$ can no longer permute the different summands. This is essentially a consequence of continuity of the semigroup. However, we
need to be  a bit cautious here for the following reason. Note that even if ${\tau _t}$ is continuous in $t$, the peripheral Poisson boundary $\mathcal{P}(\tau _t)$ is typically not continuous at zero, as $\mathcal{P}(\tau _0)$ is always the whole space $\mathbb{M}_d$,
but $P(\tau _t)$ for $t>0$ is not the whole space unless we are having automorphisms.
 
\begin{theorem}\label{deco}
Let $\tau=\{\tau_t\}_{t\in\mathbb{R}_+}:\mathbb{M}_d\to \mathbb{M}_d$ be a stationary quantum Markov semigroup. Then, there exists a decomposition of the Hilbert space $\mathcal{H}=\mathbb{C}^d$ of the form $\mathcal{H}=\bigoplus_{k=1}^M\mathcal{H}_{k,1}\otimes\mathcal{H}_{k,2}$, where $\mathcal{H}_{k,1}$, $\mathcal{H}_{k,2}$, are Hilbert spaces, such that $\mathcal{P}(\tau)=\bigoplus_{k=1}^m\mathscr{B}(\mathcal{H}_{k,1})\otimes 1_{k,2}$. Therefore, every element of $\mathcal{P}(\tau)$ is of the form $X=\bigoplus_{k=1}^m X_k \otimes 1_{k,2}$, and 
$\tau_t(X)=\bigoplus_{k=1}^m U_t^{(k)*}X_{k}U_t^{(k)} \otimes 1_{k,2}$, where $U^{(k)}$ are one parameter unitary semigroups on $\mathcal{H}_{k,1}$ for all $k$. 
\end{theorem}
\begin{proof} 
For $t>0$,  $\mathcal{P}(\tau)=\mathcal{P}(\tau_t),$ 
independent of $t$ by Theorem \ref{ppbQ}. So we appeal to the previous theorem and we have a decomposition of $\mathbb{C}^d$ not depending on time. Due to stationarity the extended Choi-Effros product $\circ_\tau$ and $\circ_{\tau_t}$, are all same as the original product.  Therefore, the Hilbert space  decomposition ($\mathcal{H}=\mathbb{C}^d$) is as mentioned, and so is the form for $\mathcal{P}(\tau)$.   
For every $t>0$, there exist unitaries $V^k_t$'s, such that  $\tau_t(X)=\bigoplus_{k=1}^mV_t^{\sigma_t(k)*}X_{\sigma_t(k)}V_t^{\sigma_t(k)} \otimes 1_{k,2}$, where for every $t>0$, $\sigma_t$ is a permutation on $\{1,2,\ldots , m\}$ such that $\text{dim}(\mathcal{H}_{k,1})=\text{dim}(\mathcal{H}_{\sigma_t(k),1})$. We show that $\sigma_t(k)=k, \ \forall t,k$. We prove this by assuming the contrary and arriving at a contradiction. Let $\sigma_{t_0}(k_0)\neq k_0, \ \mbox{for some} \ t_0>0 ~\mbox{and}~1\leq k_0\leq m$. Let $X_{k_0}=1_{k_0,1}$ and rest of the $X_k$'s be 0, which means $X_{\sigma_{t_0}(k_0)}=0$. The continuity of $\tau_t$ and the finiteness of the permutation group obviously leads to a contradiction.  Now the result follows. \end{proof}

Now we explore the non-stationary case. The discrete semigroup situation has been worked out in \cite{DA}. 
Given a unital completely positive map $\tau:\mathbb{M}_d\to\mathbb{M}_d$, now the Hilbert space $\mathcal{H}=\mathbb{C}^d$,
will decompose as  $\mathcal{H}= \mathcal{H}_0 \oplus \mathcal{H}_1$, where $\mathcal{H}_0$, itself
decomposes as $\mathcal{H}_0=
\bigoplus_{k=1}^m\mathcal{H}_{k,1}\otimes\mathcal{H}_{k,2}$ and on this space the structure is similar to the stationary case.
One of the main insights of \cite{DA} is that decomposition $\mathcal{H}=\mathcal{H}_0\oplus \mathcal{H}_1$, can be
got by studying the adjoint  map $\tau ^* ,$ which is a trace preserving completely positive map. The peripheral eigenspace of
it can be defined analogously. We have the projection of $\mathbb{M}_d$ onto the peripheral eigenspace denoted by $\mathcal{P}_P(\tau ^*).$
The space $\mathcal{H}_0$ is the projection onto the range of $\mathcal{P}_P(\tau ^*)(1)$ and $\mathcal{H}_1$ is nothing but the kernel of the matrix $\mathcal{P}_P(\tau ^*) (1).$
It has been noted in ref. \cite{DA} that it is possible to obtain a channel with a strictly positive fixed point, by restricting $\tau^*$ to $\mathscr{B}(\mathcal{H}_0)$.
Define $\gamma _{00}:\mathscr{B}(\mathcal{H}_0)\to \mathscr{B}(\mathcal{H}_0)$ as follows: for 
$X_{00}\in\mathscr{B}(\mathcal{H}_0)$, $\tau^*(\begin{bmatrix}
X_{00}&0\\0&0
\end{bmatrix})=\begin{bmatrix}
\gamma _{00}(X_{00})&0\\0&0
\end{bmatrix}$. Since $\tau^*$ is CPTP, naturally $\gamma _{00}$ is also CPTP.
Define $\tau_{00}:\mathscr{B}(\mathcal{H}_0)\to \mathscr{B}(\mathcal{H}_0)$ as $\tau_{00}=\gamma^*_{00}$, which is clearly a UCP map.
Finally, one has the following structure theorem.

\begin{theorem}(\cite{DA}, \cite{WP},  \cite{WMM})
For the UCP map $\tau:\mathbb{M}_d\to\mathbb{M}_d$, it is possible to find the decomposition of the Hilbert space $\mathcal{H}=(\bigoplus_{k=1}^m\mathcal{H}_{k,1}\otimes\mathcal{H}_{k,2})\oplus \mathcal{H}_1$, where $\mathcal{H}_1=ker(\mathcal{P}_P^*(1))$, such that there exists a CP  map $\pi_{11}:\bigoplus_{k=1}^m\mathscr{B}(\mathcal{H}_{k,1})\otimes 1_{k,2}\to \mathscr{B}(\mathcal{H}_1)$ and the peripheral Poisson boundary of $\tau$ is $$\mathcal{P}(\tau)=\{X\oplus \pi_{11}(X): \ X\in\bigoplus_{k=1}^m\mathscr{B}(\mathcal{H}_{k,1})\otimes 1_{k,2}\}.$$ The action of $\tau$ on $\mathcal{P}(\tau)$ is defined as $$\tau(X\oplus \pi_{11}(X))= \tau_{00}(X)\oplus \pi_{11}(\tau_{00}(X)),$$ for every $X \in \bigoplus_{k=1}^{m}\mathscr{B}(\mathcal{H}_{k,1})\otimes 1_{k,2}$, where $\tau_{00}$ restricted on $\bigoplus_{k=1}^{m}\mathscr{B}(\mathcal{H}_{k,1})\otimes 1_{k,2}$ is $\tau_{00}(\bigoplus_{k=1}^m X_{k} \otimes 1_{k,2})=\bigoplus_{k=1}^m U^*_{\sigma(k)}X_{\sigma(k)}U_{\sigma(k)} \otimes 1_{k,2}$, where $\sigma$ is a permutation on the $\{1,2,\ldots , m\}$, such that the $\dim (\mathcal{H}_{k,1})=\dim (\mathcal{H}_{\sigma(k),1})$.
\end{theorem}


In Theorem \ref{deco}, due to the assumption of the stationarity condition on the quantum Markov semigroup $\tau=\{\tau_t\}_{t\in\mathbb{R}_+}$, the Hilbert space $\mathcal{H}_1$ is absent. Therefore, the peripheral Poisson boundary was solely determined by the decomposition of $\mathcal{H}_0$. The following theorem strengthens Theorem \ref{deco}, as we drop the stationarity condition, with the peripheral Poisson boundary being influenced also by the Hilbert space $\mathcal{H}_1$.

\begin{theorem}\label{general}
Let $\tau=\{\tau_t\}_{t\in\mathbb{R}_+}:\mathbb{M}_d\to \mathbb{M}_d$ be a quantum Markov semigroup. Then, there exists a decomposition of the Hilbert space $\mathcal{H}=\mathbb{C}^d$ of the form $\mathcal{H}=(\bigoplus_{k=1}^m\mathcal{H}_{k,1}\otimes\mathcal{H}_{k,2})\oplus \mathcal{H}_1$, where $\mathcal{H}_1=ker(\mathcal{P}_P^*(1))$, and $\mathcal{H}_{k,i}$ for $i=1,2$, are Hilbert spaces, such that $$\mathcal{P}(\tau)=\{X\oplus \pi_{11}(X): \ X\in\bigoplus_{k=1}^m\mathscr{B}(\mathcal{H}_{k,1})\otimes 1_{k,2}\},$$ where $\pi _{11}:\mathscr{B}(\mathcal{H}_{k,1})\otimes 1_{k,2}\to \mathscr{B}(\mathcal{H}_1)$ is a CP map. The action of the semigroup $\tau$ on $\mathcal{P}(\tau)$ is 
$$\tau_t(X\oplus \pi_{11}(X))=(\bigoplus_{k=1}^m U_t^{(k)*}X_{k}U_t^{(k)} \otimes 1_{k,2})\oplus \pi_{11}(\bigoplus_{k=1}^m U_t^{(k)*}X_{k}U_t^{(k)} \otimes 1_{k,2}),$$ for every $X=\bigoplus_{k=1}^{m} X_k\otimes 1_{k,2} \in \bigoplus_{k=1}^{m}\mathscr{B}(\mathcal{H}_{k,1})\otimes 1_{k,2}$, with some unitary representations $(U^{(k)},\mathcal{H}_{k,1}), \ \forall 1\leq k\leq m$.
\end{theorem}
\begin{proof}Since, $\mathcal{ P}^*_P=\mathcal{P}_P(\tau _t^*)$, the decomposition $\mathcal{H}=\mathcal{H}_0\oplus\mathcal{H}_1$ remains intact throughout $t>0$.  Restriction of $\tau^*_t$ to $\mathscr{B}(\mathcal{H}_0)$, is a CPTP map with a strictly positive fixed point. Name the map as $\gamma ^{00}_t:\mathscr{B}(\mathcal{H}_0)\to\mathscr{B}(\mathcal{H}_0)$, for every $t>0$. 
Let $\tau^{(00)}_t=(\gamma _t^{00})^*:\mathscr{B}(\mathcal{H}_0)\to\mathscr{B}(\mathcal{H}_0)$, then, $\tau^{(00)}=\{\tau^{(00)}_t\}_{t\in\mathbb{R}_+}$ is clearly a stationary quantum Markov semigroup. Let $\pi_{00}$ and $\pi^{(t)}_{00}$ be the peripheral projection onto $\mathcal{P}(\tau^{(00)})$ and $\mathcal{P}(\tau_t^{(00)})$ respectively. Hence, $\pi_{00}=\pi^{(t)}_{00}, ~\forall t>0$. Then, due to stationarity of the quantum Markov semigroup $\tau^{(00)}_t$, $\mathcal{P}(\tau^{(00)})=\mathcal{P}(\tau_t^{(00)})$, and there exists a decomposition of the Hilbert space $\mathcal{H}_0$ of the form $\mathcal{H}_0=\bigoplus_{k=1}^m\mathcal{H}_{k,1}\otimes\mathcal{H}_{k,2}$, and $\mathcal{H}_{k,i}$ for $i=1,2$, are Hilbert spaces, such that $\mathcal{P}(\tau^{(00)})=\mathcal{P}(\tau_t^{(00)})=\bigoplus_{k=1}^m\mathscr{B}(\mathcal{H}_{k,1})\otimes 1_{k,2}$.  Therefore, by Theorem \ref{deco}, $\tau^{(00)}_t(X)=\bigoplus_{k=1}^m U_t^{(k)*}X_{k}U_t^{(k)} \otimes 1_{k,2}$ for $X=\bigoplus_{k=1}^m X_k\otimes 1_{k,2}\in\mathcal{P}(\tau^{(00)})$, with unitary representations $(U^{(k)},\mathcal{H}_{k,1}), \ \forall k$. Since $\forall t>0 ,~\mathcal{P}_P=\mathcal{P}_P(\tau _t)$, the maps $\pi _{{11}}^{(t)}$ of $\tau _t$ are independent of $t$. We may call this CP map as $\pi_{{11}}$. Hence, we are done. 
\end{proof}


We illustrate the above theory with the following example.

\begin{example}
Consider the quantum Markov semigroup $\tau=\{\tau_t\}_{t\in\mathbb{R}_+}:\mathbb{M}_3\to\mathbb{M}_3$ with a generator $\mathcal{L}$ defined on $\mathbb{M}_3$ as $\mathcal{L}(X)=E_{31}XE_{13}-\frac{1}{2}(E_{33}X-XE_{33}), ~\forall X\in\mathbb{M}_3$, where $E_{ij}$'s are matrix units.  That is, 
\begin{eqnarray*}
    \mathcal{L}(\begin{bmatrix}
        x_{11}&x_{12}&x_{13}\\
        x_{21}&x_{22}&x_{23}\\
        x_{31}&x_{32}&x_{33}
    \end{bmatrix})=\begin{bmatrix}
        0&0&-\frac{1}{2}x_{13}\\
        0&0&-\frac{1}{2}x_{23}\\
        -\frac{1}{2}x_{31}&-\frac{1}{2}x_{32}&x_{11}-x_{33}
    \end{bmatrix}.
\end{eqnarray*}
Let $\tau^*=\{\tau_t^*\}_{t\in\mathbb{R}_+}$ be the semigroup of CPTP maps, where each $\tau^*_t$ is the Hilbert-Schmidt adjoint of $\tau_t$ for every
$t>0$. Clearly, $\tau^*_t=\exp{(t\mathcal{L^*})}, ~\forall t>0$, where $\mathcal{L}^*$ is the adjoint of the generator $\mathcal{L}$.  Direct computations show that \begin{eqnarray*}
    \mathcal{L}^*(\begin{bmatrix}
        x_{11}&x_{12}&x_{13}\\
        x_{21}&x_{22}&x_{23}\\
        x_{31}&x_{32}&x_{33}
    \end{bmatrix})=\begin{bmatrix}
        x_{33}&0&-\frac{1}{2}x_{13}\\
        0&0&-\frac{1}{2}x_{23}\\
        -\frac{1}{2}x_{31}&-\frac{1}{2}x_{32}&-x_{33}
    \end{bmatrix}.
\end{eqnarray*}
The peripheral Poisson boundary for $\tau^*$ is $\mathcal{P}(\tau^*)=\mathbb{M}_2\oplus 0$. Therefore, $\mathcal{H}_0=span\{e_1,e_2\}$ and $\mathcal{H}_1=\mathbb{C}e_3$. Turns out the quantum Markov semigroup $\tau^{(00)}=\{\tau_t^{(00)}\}_{t\in\mathbb{R}_+}:\mathbb{M}_2 \to\mathbb{M}_2$ is $\tau^{(00)}_t=id., ~\forall t>0$. Clearly, $\mathcal{P}(\tau^{(00)})=\mathbb{M}_2$ and the CP map $
\pi_{11}:\mathbb{M}_2\to \mathbb{C}$ is $$\pi_{11}(\begin{bmatrix}
x_{11}&x_{12}\\
x_{21}&x_{22}
\end{bmatrix})=x_{11}.$$
The peripheral Poisson boundary of $\tau$ is $$\mathcal{P}(\tau)=\{\begin{bmatrix}
    x_{11}&x_{12}&0\\
        x_{21}&x_{22}&0\\
        0&0&x_{11}
\end{bmatrix}, ~ x_{11},x_{12},x_{21},x_{22}\in\mathbb{C}\}.$$ Moreover, $\tau$ restricted to $\mathcal{P}(\tau)$ is just the identity map.
\end{example}

\begin{remark} 
Let $\tau:\mathbb{M}_d\to \mathbb{M}_d$ be a unital trace preserving completely positive map. From ref. \cite{DA}, it is known that there exists a decomposition of the  Hilbert space $\mathcal{H}=\mathbb{C}^d$ of the form $\mathcal{H}=\bigoplus_{k=1}^m \mathcal{H}_{k,1}\otimes \mathcal{H}_{k,2}$, where  $\mathcal{H}_{k,1}$, $\mathcal{H}_{k,2}$  are Hilbert spaces, and the peripheral Poisson boundary of $\tau^*$ is given as $\mathcal{P}(\tau^*)=\bigoplus_{k=1}^m \mathscr{B}(\mathcal{H}_{k,1})\otimes \rho_k$, where  $\rho_k$  is some invertible state in $\mathscr{B}(\mathcal{H}_{k,2})$. The peripheral projection of $\tau$ is  $$\mathcal{P}_P(X)=\bigoplus_{k=1}^m tr_{k,2}[(P_kXP_k)(1_{k,1}\otimes\rho_k)]\otimes 1_{k,2},$$ where $tr_{k,2}$ denotes the trace in the $\mathcal{H}_{k,2}$ factor.  Being a limit of trace preserving maps, $\mathcal{P}_P$ is also trace preserving. Then  it follows that $\rho _k=\frac{1}{\text{dim}{(\mathcal{H}_{k,2}})}1_{\mathcal{H}_{k,2}}.$ Analogous results hold  even in the continuous time, that is, for $\mathcal{P}(\tau _t)$ with $t>0$, for trace preserving quantum Markov semigroups $\{\tau _t: t\in \mathbb{R}_+\}$.

\end{remark}

\section*{Acknowledgments}
Bhat gratefully acknowledges funding from  ANRF  (India) through J. C.
Bose Fellowship No. JBR /2021/000024.
Dias is thankful to Indian Statistical Institute, Bangalore for the kind hospitality during her visits.

\textbf{AI declaration:} The AI (Gemini Pro) was used to verify and make minor corrections of computations in Example \ref{unitary mix}
and Example \ref{Grigorchuk}.


\begin{thebibliography}\bibliographystyle
\bibitem{WA} W. B. Arveson, Sub-algebras of $C^*$-algebras, Acta Math. 123 (1969), 141-224.

\bibitem{DA} D. Amato, P. Facchi, A. Konderak, Asymptotic Dynamics in the Heisenberg Picture: Attractor Subspace and Choi-Effros Product (2025).

\bibitem{BKT1}  B. V. R. Bhat, S. Kar and B. Talwar, Peripheral
Poisson boundary,  Isr. J. Math. (2025).  https://doi.org/10.1007/s11856-025-2830-2


\bibitem{BKT2} B. V. R. Bhat, S. Kar and B. Talwar, Peripherally
automorphic unital completely positive maps, Linear Alg. and its
Appli. {\bf 678} (2023)191-205.

\bibitem{BA1} B. V. R. Bhat, A. S. Dias, Peripheral Poisson boundary: Extensions and examples, Proc. Math. Sci. {\bf 136}, 28 (2026).  


\bibitem{BhC} B. V. R. Bhat, Minimal dilations of quantum dynamical semigroups to semigroups of
endomorphisms of $C^*$-algebras, Journal of the Ramanujan Mathematical Society, {\bf 14},
(1999), 109–124.

\bibitem{BVRB} B. V. R. Bhat, Cocycles of CCR Flows, Memoirs of American Mathematical Society, No. 709  (2001).

\bibitem{NBNO} N. Brown, N. Ozawa, $C^*$-Algebras and Finite Dimensional Approximations, Vol. 88, (2008).

\bibitem{BGKS} A. Batkai, U. Groh, D. Kunszenti-Kovacs, M. Schreiber, Decomposition of operator semigroups on $W^*$-algebras, DOI 10.1007/s00233-011-9362-2, Springer, (2011).


\bibitem{Emel} Eduard Yu. Emel'yanov, Non-spectral Asymptotic Analysis of One-Parameter Operator Semigroups, Chapter 1, (2007).

\bibitem{EFHN} T. Eisner, B. Farkas, M. Haase, R. Nagel, Operator Theoretic Aspects of Ergodic Theory, vol. 272, Springer-Verlag (2015).


\bibitem{RGIP} R. Grigorchuk, I. Pak, Groups of Intermediate Growth, an introduction, (2008).


\bibitem{Izumi2002} M. Izumi, Non-commutative Poisson boundaries and compact quantum group actions, Adv. Math. {\bf 169} (2002), no.~1, 1--57. MR1916370

\bibitem{Izumi2004} M. Izumi, Non-commutative Poisson boundaries, Discrete Geometric Analysis, 69-81, Contemp. Math. 347 Amer. Math. Soc. Providence, RI 2004.


\bibitem{Lind} G. Lindblad, A general no-cloning theorem, Letters to Mathematical Physics 47: 189-196, (1999).

\bibitem{BP} B. Prunaru, A short proof of Choquet-Deny theorem for abelian topological semigroups (2000).

\bibitem{VP} V. Paulsen, Completely Bounded Maps and Operator Algebras, Chapter 4, Cambridge Studies in Advanced Mathematics (2003).

\bibitem{CS} C. R. Rao, D. N. Shanbhag, An elementary proof for an extended version of the Choquet-Deny theorem, Journal of Multivariate Analysis 38, 141-148, (1991).

\bibitem{Watrous}  J. Watrous, The Theory of Quantum Information, Cambridge University Press, (2018).

\bibitem{WP} M. M. Wolf and D. Perez-Garcia, The inverse eigenvalue problem for quantum channels
(arXiv:1005.4545 [quant-ph]) (2010).

\bibitem{WMM} M. M. Wolf, Quantum Channels and Operations: A Guided Tour. Lecture Notes, TUM, (2012).

\end{thebibliography}
\end{document}